\documentclass[11pt]{article}

\usepackage{amsmath,amssymb,amsthm,setspace,graphicx,array,tikz,url,blkarray,easyReview,bm}
\usepackage[text={6.5in,8.4in}, top=1.3in]{geometry}
\newtheorem{theorem}{Theorem}[section]
\newtheorem{proposition}[theorem]{Proposition}
\newtheorem{lemma}[theorem]{Lemma}
\newtheorem{corollary}[theorem]{Corollary}

\newtheorem{problem}[theorem]{Problem}
\newtheorem{definition}[theorem]{Definition}
\newtheorem{example}[theorem]{Example}

\newcommand{\G}{{\mathcal G}}
\renewcommand{\P}{{\mathcal P}}

\begin{document}
	
\title{Inducibility in $H$-free graphs and inducibility of Tur\'an graphs}
	
\author{
	Raphael Yuster
	\thanks{Department of Mathematics, University of Haifa, Haifa 3103301, Israel. Email: raphy@math.haifa.ac.il\;.}
}

\date{}
	
\maketitle
	
\setcounter{page}{1}
	
\begin{abstract}
	For graphs $F$ and $H$, let $i(F)$ denote the inducibility of $F$ and let $i_H(F)$ denote the inducibility of $F$ over $H$-free graphs.
	We prove that for almost all graphs $F$ on a given number of vertices, $i_{K_k}(F)$ attains infinitely many values as $k$ varies.
	For complete partite graphs $F$ (and, more generally, for symmetrizable families of graphs $F$), we prove that 
	$i_H(F)=i_{K_k}(F)$ where $k=\chi(H)$, and is attained by a complete $\ell$-partite graphon $W_{F,k}$, where $\ell < k$.
	
	We determine the part sizes of $W_{F,k}$ for
	all $k$, whence determine $i(F)$, whenever $F$ is the Tur\'an graph on $s$ vertices and $r$ parts, for all $s \le 3r+1$, which was recently proved by Liu, Mubayi, and Reiher for $s=r+1$. As a corollary, this determines the inducibility of all Tur\'an graphs on at most $14$ vertices.
	Furthermore, since inducibility is invariant under complement, this determines the inducibility of all
	matchings and, more generally, all graphs with maximum degree $1$, of any size.
	Similarly, this determines the inducibility of all triangle factors, of any size.
	
	For complete partite graphs $F$ with at most one singleton part, we prove that
	$i_{K_k}(F)$ only attains finitely many values as $k$ varies; in particular, there exists $t=t(F)$ such that $i(F)$ is attained by some complete $t$-partite graphon. This is best possible as it was shown by 
	Liu, Pikhurko, Sharifzadeh, and Staden that this is not necessarily true if there are two singleton parts.
	
	Finally, for every $r$, we give a nontrivial sufficient condition for a complete $r$-partite graph $F$ to have the property that $i(F)$ is attained by a complete partite graphon all whose part sizes are distinct.
	

\end{abstract}

\section{Introduction}\label{sec:intro}

In the generalized inducibility problem we are given two finite sets of graphs: a {\em forbidden} set $X$, and a {\em target} set $F$ where all elements of $F$ have the same order. The goal is to determine the maximum number of induced copies of elements of $F$ in a large $X$-free\footnote{Recall that a graph is $H$-free if it does not contain $H$ as a (not necessarily induced) subgraph.} graph. If $F$ is closed under edge addition, this is the generalized Tur\'an problem \cite{AS-2016}.
When $(X,F)=(\{K_k\},\{K_2\})$ this is Tur\'an's Theorem \cite{turan-1941}
and when $X = \emptyset$, this is the well-known inducibility problem \cite{PG-1975}.

More formally, for graphs $F$ and $G$, let $P(F,G)$ be the number of induced copies of $F$ in $G$, namely the number of $|V(F)|$-subsets of $V(G)$ that induce a subgraph isomorphic to $F$.
Let $p(F,G)=P(F,G)/\binom{|V(G)|}{|V(F)|}$ be the {\em induced density} of $F$ in $G$.
These definitions naturally extend to the case where $F$ is a family of graphs of the same order.
Let $i(F,n)$ be the maximum of $p(F,G)$ taken over all graphs $G$ with $n$ vertices and let
$i_X(F,n)$ be the maximum of $p(F,G)$ taken over all $X$-free graphs $G$ with $n$ vertices.
When $X=\{H\}$ is a singleton, we use the notation $i_H(F,n)$.
Finally, let $i(F) = \lim_{n \rightarrow \infty} i(F,n)$ and
let $i_H(F) = \lim_{n \rightarrow \infty}  i_H(F,n)$. It is easily seen by monotonicity that these limits exist and are equal to the maximum density of induced copies of $F$ in a graphon. We may sometimes focus on the special important case $H=K_k$ so we denote $i_k(F,n)=i_{K_k}(F,n)$ and $i_k(F)=i_{K_k}(F)$.
The parameter $i(F)$ is the {\em inducibility} of $F$ and the parameter $i_H(F)$ is, therefore, the inducibility of $F$ over $H$-free graphs.
Notice that $i(F)=i(F^c)$ where $F^c$ is the complement of $F$, but this is not necessarily so for $i_H(F)$.

Generalized inducibility problems date back to results of Zykov \cite{zykov-1949}
and Erd\H{o}s \cite{erdos-1962} who resolved $i_k(K_r)$. In fact, their proofs show that $i_k(K_r,n)$ is attained by the Tur\'an graph $T(n,k-1)$ (notice that this also holds in the trivial case where $k \le r$). The systematic study of generalized Tur\'an problems was
initiated by Alon and Shikhelman \cite{AS-2016}. Inducibility (of non-complete graphs) dates back
to a classical result of Goodman \cite{goodman-1959} who proved $i(K_{1,2})=\frac{3}{4}$
but the systematic study of the parameter $i(F)$ was initiated by
Pippenger and Golumbic \cite{PG-1975}.
Since then, many strong results concerning inducibility as well as inducibility in $H$-free graphs have been proved, but the problem at large is still an unclaimed land. We briefly mention some of the main results in this area and later expand on those that are particularly relevant to our results.

\subsection{Related research}

Exoo \cite{exoo-1986} determined the inducibilities of some graphs on four vertices and
provided bounds for all of them. Yet, even  the inducibility of $P_4$, the path on four vertices, is still undetermined, see \cite{EL-2015,flagmatic-site} for the best known lower and upper bounds.
The inducibilities of the remaining graphs on four vertices were determined by Hirst \cite{hirst-2014}.
The inducibility of complete bipartite graphs was determined by Brown and Sidorenko \cite{BS-1994};
in fact, they have proved that $i(K_{s,t},n)$ is attained by a complete bipartite graph.
Bollob\'as, Egawa, Harris, and Jin \cite{BEHJ-1995} determined the inducibility of complete $r$-equipartite graphs $F$ assuming that the part sizes are not too small with respect to $r$;
in fact, they proved that under such conditions, $i(F,n)$ is uniquely attained by the $r$-partite Tur\'an graph with $n$ vertices, for $n$ sufficiently large, strengthening another result from \cite{BS-1994}. Hatami, Hirst, and Norine proved that the inducibility of graphs that are large blow-ups is also attained by blow-ups \cite{HHN-2014}. In their paper, Pippenger and Golumbic \cite{PG-1975}  observed the generic lower bound $i(F) \ge r!/(r^r-r)$ where $r=|V(F)|$, and conjectured that this is tight
for all cycles on at least five vertices. This was proved many years later for $C_5$ by Balogh, Hu, Lidick{\`y}, and Pfender \cite{BHLP-2016}. Some progress on the general conjecture has been obtained in \cite{HT-2018,KNV-2019}. It was proved by Fox, Huang, and Lee \cite{FHL-2017} and by the author \cite{yuster-2019} that almost all graphs have $i(F) = r!/(r^r-r)$. Recently, Ueltzen \cite{ueltzen-2024} classified all graphs with high inducibility. Stability results for the induced density of certain graphs, as well as the inducibility of some graphs on five vertices were obtained by Pikhurko, Slia{\v{c}}an, and Tyros \cite{PST-2019}. Further stability results, as well as the inducibility of $K_{2,1,1,1}$ and $K_{3,1,1}$ were obtained by Liu, Pikhurko, Sharifzadeh, and Staden \cite{LPSS-2023}.
Recently, Liu, Mubayi, and Reiher \cite{LMR-2023} determined the inducibility of $K_r^{-}$ for all $r$.

Inducibility of graph families is also an intriguing area of study.
Most notably, we mention the {\em edge statistics} problem, asking for the limiting behavior of inducibility of the family of $r$-vertex graphs having $\ell$ edges, as $\ell$ and $r$ vary. The study of edge statistics was initiated by Alon, Hefetz, Krivelevich, and Tyomkyn \cite{AHKT-2020} and the problem was completely resolved by the combined results of
Kwan, Sudakov, and Tran \cite{KST-2019}, Martinsson, Mousset, Noever, and Truji{\'c} \cite{MMNT-2019} and
Fox and Sauermann \cite{FS-2020}.

Finally we mention some more recent results on inducibility over $H$-free graphs.
Answering a question of Erd\H{o}s, Grzesik \cite{grzesik-2012} and
Hatami, Hladk{\`y}, Kr{\'a}l, Norine, and Razborov \cite{HHKNR-2013} determined $i_3(C_5)$. A strong form of this result, determining $i_3(C_5,n)$ for every $n$, was obtained by Lidick{\`y} and Pfender \cite{LP-2018}. Another question of Erd\H{o}s asks for
$i_k(F_r,n)$ and $i_k(F_r)$ where $F_r$ is the set of all graphs on $r$ vertices excluding the empty graph. The case $(k,r)=(3,3)$ follows from Goodman's result \cite{goodman-1959},
the case $k=3$ and $r=4,5$ was determined by Das, Huang, Ma, Naves, and Sudakov \cite{DHMNS-2013} and the case $k=3$ and $r=4,5,6,7$ was determined by Pikhurko and Vaughan \cite{PV-2013}.

\subsection{Our results}

Our aim is to initiate a more systematic study of inducibility over $H$-free graphs. As we shall see, this approach is beneficial since in some cases it streamlines to obtain new results on inducibility proper.

Our first main result proves that $i_H(F)$, and in fact already $i_k(F)$, is usually ``restricting''.
By this we mean that for almost all graphs $F$ (in a well-defined sense made precise), $i_k(F)$ attains infinitely many values as $k$ varies.

\begin{definition}\label{def:diverse}
Call a graph $F$ {\em inducibility diverse} if for all $k \in \mathbb{N}$, there exists $m > k$ such that
$i_k(F) < i_m(F)$. In particular, $i_k(F) < i(F)$ for every $k$.
\end{definition}
One may wonder why it is not obvious that most $F$ are inducibility diverse? It is quite easy to see that complete graphs $F$ are inducibility diverse, as $i(F,n)=1$ is only attained by $K_n$. However, this is already false for $F$ that are complete graphs missing a single edge --- they are {\em not} inducibility diverse, as proved by Liu, Mubayi, and Reiher  \cite{LMR-2023}. In fact, for every non self-complementary four-vertex graph,
either it or its complement is not inducibility diverse, see Even-Zohar and Linial \cite{EL-2015}.

Recall that a random graph on $r$ vertices is the probability distribution $\G(r,\frac{1}{2})$ on $r$-vertex graphs where each pair of vertices is an edge with probability $\frac{1}{2}$ and all $\binom{r}{2}$ choices are independent. We say that $F \sim \G(r,\frac12)$ satisfies a graph property {\em asymptotically almost surely}
if $F$ satisfies the property with probability $1-o_r(1)$.

\begin{theorem}\label{t:fractil}
	Asymptotically almost surely, it holds that $F \sim \G(r,\frac12)$ is inducibility diverse.
\end{theorem}

Our remaining results are about $i_H(F)$ where $F$ is a complete partite graph or, more generally, $F$ is a family of symmetrizable graphs (see definition below). Inducibility of complete partite graphs has been widely studied with many strong results, yet still unresolved in most cases.
This is in spite of the relatively simple fact that $i(F,n)$ ``just amounts'' to solving a polynomial optimization problem. These graphs are further particularly appealing in our case since for such $F$ we can show that $i_H(F)=i_k(F)$ where $k=\chi(H)$, see Theorem \ref{t:symmetrizable} below.

To state our results, we need several definitions.
Recall that a graphon is a symmetric function $W: [0,1]^2 \rightarrow \{0,1\}$.
Any set $V$ of distinct points in $[0,1]$ corresponds to an induced graph on $V$ in the graphon $W$:
for distinct $u,v \in V$ we have that $uv$ is an edge if and only if $W(u,v)=1$.
For a graph $H$ and a graphon $W$, the {\em induced density of $H$ in $W$} is the probability that 
a randomly chosen set of $|V(H)|$ points in $[0,1]$ induces a copy of $H$ in $W$.
For more information on graphons and graph limits, see Lov\'asz \cite{lovasz-2012}.
Let $\Delta^{\ell-1}$ denote the standard $(\ell-1)$-dimensional simplex and
let $(x_1,\ldots,x_\ell) \in \Delta^{\ell-1}$. Partition $[0,1]$ into $\ell$ intervals of sizes
$x_1,\ldots,x_\ell$ and let $W=W[x_1,\ldots,x_\ell]: [0,1]^2 \rightarrow \{0,1\}$ be the graphon defined as follows:
$W(x,y)=1$ if and only if $x$ and $y$ are in distinct intervals.
Most importantly, 
notice that if $\chi(H) > \ell$, then $W$ is $H$-free.
It is easily seen that $W$ is the limit of a sequence of complete $\ell$-partite graphs,
where the $n$th element in the sequence has $n$ vertices and the size of its $i$th part is either $\lfloor x_i n \rfloor$ or $\lceil x_i n \rceil$ (we shall sometimes allow
some of the $x_i$'s to be zero, so if there are only $\ell^*$ nonzero entries in $(x_1,\ldots,x_\ell)$, then the sequence elements are, in fact, complete $\ell^*$-partite graphs).
We call $W[x_1,\ldots,x_\ell]$ a {\em complete $\ell$-partite graphon}.

We recall the notion of {\em Zykov symmetrization operation} (hereafter {\em symmetrization}) \cite{zykov-1949}.
Suppose that $x,y$ are two nonadjacent vertices of a graph $F$. Obtain a graph $F'$ by replacing $x$ with a copy of $y$, keeping $y$ and its clone nonadjacent (so, following a symmetrization,
$y$ and its clone are {\em twins}). We say that a family of graphs is {\em symmetrizable} if it is closed 
under symmetrization. Clearly, every family of complete partite graphs is symmetrizable, and clearly the union and intersection of symmetrizable families is symmetrizable. But note that there are other nontrivial symmetrizable families of graphs. For example, $\{K_{1,3}, K_4^-,D\}$ is symmetrizable, where $D$ is the (non complete partite) graph on four vertices obtained by adding an edge to $K_{1,3}$ (sometimes $D$ is called the {\em paw graph}). It is easy to see from the definition that complete partite graphs are just the singleton symmetrizable families.

The following theorem makes $i_H(F)$ more tractable to determine, whenever $F$ is a family of symmetrizable graphs.
\begin{theorem}\label{t:symmetrizable}
	Let $F$ be a symmetrizable family of graphs. Then for any graph $H$ it holds that $i_H(F)=i_{\chi(H)}(F)$.
	In particular, $i_H(F)$ is a attained by a complete $\ell$-partite graphon where $\ell < \chi(H)$.
\end{theorem}
\noindent 
Theorem \ref{t:symmetrizable} is, in general, false for non-complete partite graphs.
For example, consider the case $F=C_5$. Taking a blowup of $C_5$ (which is triangle free) already gives that $i_3(C_5) \ge 24/625$ (in fact, this is tight \cite{grzesik-2012,HHKNR-2013}) but it is not difficult to show that, say, $i_{C_7}(C_5)=0$, yet $\chi(C_7)=3$.

For a symmetrizable family of graphs $F$, and for a graph $H$ with $\chi(H)=k$, Theorem \ref{t:symmetrizable} asserts that in order to determine $i_H(F)$, it suffices to determine a point
$(x_1,\ldots,x_\ell) \in \Delta^{\ell-1}$ where $\ell < k$ such that $W[x_1,\ldots,x_\ell]$ is a complete $\ell$-partite graphon attaining $i_k(F)$. Denote such a graphon by $W_{F,k}$ \footnote{Although the point $(x_1,\ldots,x_\ell)$
is not necessarily unique, it is immediate that the induced density of $F$ in $W_{F,k}$, namely $i_k(F)$, is determined by the point, so it will not be confusing to use the notation $W_{F,k}$ despite its possible non-uniqueness.}.

As it turns out, different types of complete partite graphs $F$ (which, recall, are the singleton symmetrizable families) may behave quite differently; some $F$ are inducibility diverse, while others are not. For some $F$, all part sizes of $W_{F,k}$ are equal (even if the part sizes of $F$ are not), while for others, they are all unequal.
We demonstrate these phenomena and determine $W_{F,k}$ for all $k$ (whence also the inducibility proper) in several notable cases.
It will be convenient to use the notation $K_{a_1,\ldots,a_r}$ to denote the complete $r$-partite graph with part sizes $a_1 \ge a_2 \ge \cdots a_r \ge 1$. When all parts sizes differ by at most one, we obtain
a Tur\'an graph. As usual, we denote by $T(s,r)$ the Tur\'an
graph with $s$ vertices and $r$ parts. Notice that complete equipartite graphs are a special case of Tur\'an graphs. To avoid trivial cases, we shall always assume that $k \ge 3$ (since $i_2(F)=0$ unless $F$ is the empty graph) and that $r \ge 2$ (i.e., $F$ is not the empty graph) as for the empty graph we have that $i_k(F)=1$ which is obtained by $W[1]$ (the empty graphon).

The first natural case to consider is that of complete graphs $F=K_r$.
Recall the aforementioned results of Erd\H{o}s and Zykov who proved that $i_k(K_r,n)$ is attained by the Tur\'an graph $T(n,k-1)$. We therefore obtain the following corollary.
\begin{corollary}[\cite{erdos-1962}]\label{coro:erdos}
	$W_{K_r,k}=W[\frac{1}{k-1},\ldots,\frac{1}{k-1}]$.
\end{corollary}

Naturally, the next pressing case are complete graphs missing an edge, $F=K_r^-=K_{2,1,\ldots,1}=T(r,r-1)$.
Let us recap the state of the art regarding $i(K_r^-)$. It follows from the aforementioned result of Goodman \cite{goodman-1959} that
$i(K_3^-)=\frac{3}{4}$ and that the corresponding graphon is $W[\frac{1}{2},\frac{1}{2}]$.
The case of $i(K_4^{-})=72/125$ was determined by Hirst \cite{hirst-2014} using the flag algebra method;
the corresponding graphon is $W[\frac{1}{5},\frac{1}{5},\frac{1}{5},\frac{1}{5},\frac{1}{5}]$. It was recently determined
by Liu, Pikhurko, Sharifzadeh, and Staden \cite{LPSS-2023} that $i(K_5^-)$ is attained by $W[\frac{1}{8},\ldots,\frac{1}{8}]$. Finally, the problem was solved completely by Liu, Mubayi, and Reiher \cite{LMR-2023} who proved that
$i(K_r^-)$ is attained by $W[\frac{1}{t},\ldots,\frac{1}{t}]$ where $t=\lceil (r-2)(3r+1)/6 \rceil$.
In particular, notice that this result implies that $K_r^{-}$ (unlike $K_r$) is not inducibility diverse.
Our next result establishes $i_k(K_r^-)$ for all pairs $(r,k)$. In fact, we will not need to prove it separately, as it follows as a corollary of a significantly more general result (Theorem \ref{t:turan} below).
\begin{theorem}\label{t:almost-complete}
	Let $(r,k) \ge (3,3)$ be a pair of integers and let $t=\lceil (r-2)(3r+1)/6 \rceil$.
	For all $k \le t+1$, $i_k(K_r^-)$ is attained by $W[\frac{1}{k-1},\ldots,\frac{1}{k-1}]$. For $k > t+1$,
	$i_k(K_r^-)$ is attained by $W[\frac{1}{t},\ldots,\frac{1}{t}]$.
\end{theorem}

Having a complete picture of the cases of complete and almost complete graphs, we next consider other natural classes of complete partite graphs. As already mentioned, the inducibility of complete bipartite graphs has been determined by Brown and Sidorenko \cite{BS-1994} (see also \cite{LPSS-2023} Theorem 1.6). They have proved that $i(K_{a,b})$ (assuming $ab > 1$) is attained by the bipartite graphon $W[\alpha,1-\alpha]$ where $\alpha = \arg \max_{x \in [0,1]} x^a(1-x)^b+x^b(1-x)^a$. Notice that this completely determines $i_k(F)$ for all complete bipartite graphs $F$,
as we have $i_k(F)=i(F)$ for all $k \ge 3$. In particular, all complete bipartite graphs (other than $K_2$)
are not inducibility diverse. We note that $i(K_{1,5})=\frac{2}{9}(5\sqrt{10}-14)$ is irrational, while
$i(K_{1,4})=\frac{5}{12}$ yet the corresponding $\alpha$ is $\frac{1}{2}-\frac{1}{\sqrt{12}}$.

Proceeding to Tur\'an graphs, the situation becomes more involved. Brown and Sidorenko \cite{BS-1994} (see also \cite{BEHJ-1995}) proved that if
\begin{equation}\label{e:BS}
\left(1+\frac{1}{r}\right)^s\left(1-\frac{s}{\lfloor s/r \rfloor(r+1)}\right) > 1
\end{equation}
then $i(T(s,r))$ is obtained by $W[\frac{1}{r},\ldots,\frac{1}{r}]$. Again, notice that this completely determines $i_k(T(s,r))$
for such pairs $(s,r)$ as we have in this case that $i_k(T(s,r))=i(T(s,r))$ for $k \ge r+1$ and clearly $i_k(T(s,r))=0$ for $k \le r$.
In particular, if $s/r$ is an integer, it is proved in \cite{BS-1994} that \eqref{e:BS} is also {\em necessary}.
So, for example, $i(T(12,6))$ is {\em not obtained} by a complete $6$-partite graphon.
We have also already seen that for almost complete graphs (which correspond to $T(r+1,r)$), the
situation is quite different and is not covered by \eqref{e:BS}.
Our next result significantly extends the state of the art regarding Tur\'an graphs and complete equipartite graphs, both for $i_k(T(s,r))$ and $i(T(s,r))$.
\begin{theorem}\label{t:turan}
	Let $2 \le r < s \le 3r+1$ and let $k > r$. Let $t \ge r$ be the largest integer for which
	$$
	\frac{t}{t-r}\left(1-\frac{1}{t}\right)^{s} > 1\;.
	$$
	(possibly $t=r$) and let $\ell= \min\{k-1,t\}$.
	Then, $i_k(T(s,r))$ is uniquely attained by $W[\frac{1}{\ell},\ldots,\frac{1}{\ell}]$.
	In particular, $i(T(s,r))$ is uniquely attained by $W[\frac{1}{t},\ldots,\frac{1}{t}]$.
	Furthermore, setting $s=pr+q$ where $0 \le q < r$ we have
	$$
	i_k(T(s,r)) = \frac{\ell!s!}{(\ell-r)!(r-q)!q!(p!)^r(p+1)^q\ell^{s}}
	$$
	and
	$$
	i(T(s,r)) = \frac{t!s!}{(t-r)!(r-q)!q!(p!)^r(p+1)^q t^{s}}\;.
	$$
\end{theorem}
Notice that Theorem \ref{t:turan} implies Theorem \ref{t:almost-complete} as the latter is the case $s=r+1$.
Furthermore, since $i(F)=i(F^c)$, Theorem \ref{t:turan} determines the inducibility of all
matchings and, more generally, all graphs with maximum degree $1$, of any size.
Similarly, it determines the inducibility of all triangle factors, of any size.

Another important consequence of Theorem \ref{t:turan} is that it settles the inducibility problem for
all Tur\'an graphs on at most $14$ vertices:
\begin{theorem}\label{t:le14}
	Let $F$ be a Tur\'an graph on at most $14$ vertices. Then $i(F)$ and as well as $i_k(F)$
	and $W_{F,k}$ for every $k$ are completely determined from Table \ref{table:le14}.
\end{theorem}

It seems highly interesting to determine whether the bound $s \le 3r+1$ in Theorem \ref{t:turan} can be lifted.
\begin{problem}
	Is it true that for all $2 \le r < s$, there exists $t=t(s,r)$ such that the following holds.
	For all $k \le t+1$, $i_k(T(s,r))$ is attained by $W[\frac{1}{k-1},\ldots,\frac{1}{k-1}]$ and for all $k > t+1$,
	$i_k(T(s,r))$, and hence also $i(T(s,r))$, is attained by $W[\frac{1}{t},\ldots,\frac{1}{t}]$.
\end{problem}
The following proposition shows that the smallest nontrivial case, i.e., the case $k=r+1$ 
in the aforementioned problem, does hold; namely $i_{r+1}(T(s,r))$ is attained by
$W[\frac{1}{r},\ldots,\frac{1}{r}]$.
\begin{proposition}\label{prop:turan}
	Let $2 \le r \le s$. Then $i_{r+1}(T(s,r))$ is attained by $W[\frac{1}{r},\ldots,\frac{1}{r}]$.
\end{proposition}

Our final two theorems apply to complete partite graphs in general. The first gives a necessary condition
for a complete partite graph to be inducibility diverse.
\begin{theorem}\label{t:not-diverse}
	Let $F=K_{a_1,\ldots,a_r}$. If $a_{r-1} > 1$, then there exists $t=t(F)$ such that for all $k > t$, $i_k(F)$,
	and hence also $i(F)$, is attained by a complete $t$-partite graphon. In particular, $F$ is not inducibility diverse.
\end{theorem} 
Theorem \ref{t:not-diverse} is best possible in the sense that it is generally false if two vertex classes are singletons.
For example, it follows from Theorem 1.9 in \cite{LPSS-2023} that $i(K_{3,1,1})$ is not attained by any
complete $k$-partite graphon for any finite $k$. It may be interesting to characterize all sequences
$(a_1,\ldots,a_{r-2},1,1)$ for which $F=K_{a_1,\ldots,a_{r-2},1,1}$ is not inducibility diverse.
Nevertheless, we do obtain the following immediate consequence of Theorems \ref{t:turan} and \ref{t:not-diverse}.
\begin{corollary}\label{coro:turan}
	A Tur\'an graph is inducibility diverse if and only if it is complete.
\end{corollary}
\noindent
Indeed, complete graphs are inducibility diverse, yet if $T(s,r)$ is such that $r < s \le 3r+1$
then Theorem \ref{t:turan} applies, or otherwise every part has at least three vertices, so Theorem \ref{t:not-diverse} applies.

The next theorem provides some additional information as to the structure of $W_{F,k}$. Call a complete partite graph $F=K_{a_1,\ldots,a_r}$ {\em strongly unbalanced} if for all $1 \le i < j \le r$ it holds that $(a_i-a_j)^2 > a_i + a_j$. For example: $K_{8,4,1}$ is strongly unbalanced.
Note that if $F$ is strongly unbalanced, then, by Theorem \ref{t:not-diverse}, $i_k(F)$ and also $i(F)$ are attained
by some complete partite graphon, yet the number of parts of the latter, while finite, may be {\em  larger} than $r$. Yet, interestingly, more can be said about the graphon's part sizes.
\begin{theorem}\label{t:strongly-unbalanced}
	Let $F$ be a strongly unbalanced complete partite graph. Then for all $k > r$, all part sizes of $W_{k,F}$ are distinct.
	Similarly, all part sizes of $W_F$, the complete partite graphon attaining $i(F)$, are distinct.
\end{theorem}

Finally, recall from the discussion above that there are irrational inducibilities, e.g. $i(K_{1,5})$ is such. But can it be the case that $i(F)$ is rational, yet $i_k(F)$ is irrational for some $k$?
We show that the answer is positive and determine the smallest complete partite graph for which this holds.
\begin{proposition}\label{prop:smallest}
	There are complete partite graphs $F$ for which $i(F)$ is rational yet $i_k(F)$ is irrational. The smallest such graph is $K_{3,1,1}$ for which $i(K_{3,1,1})=\frac{216}{625}$ yet $i_5(K_{3,1,1})$ is irrational.
\end{proposition}
\noindent We note that $i(K_{3,1,1})=\frac{216}{625}$ was determined recently in \cite{LPSS-2023}.

The proof sections of this paper are organized as follows. In Section \ref{sec:almost-all} we prove
Theorem \ref{t:fractil} showing that almost all graphs are inducibility diverse. In Section \ref{sec:symmetrizable} we prove Theorem \ref{t:symmetrizable} showing that $i_H(F)=i_{\chi(H)}(F)$
for symmetrizable families $F$. In section  \ref{sec:turan} we consider Tur\'an graphs and prove Theorems \ref{t:turan}, \ref{t:le14} and Proposition \ref{prop:turan}.
In Section \ref{sec:more} we consider general complete partite graphs and prove Theorems \ref{t:not-diverse},
\ref{t:strongly-unbalanced} and Proposition \ref{prop:smallest}.

\section{Almost all graphs are inducibility diverse}\label{sec:almost-all}

In this section we prove Theorem \ref{t:fractil}.
Throughout this section, $F$ is a graph and $r=|V(F)|$. We shall consider the number of copies of $F$
in graphs $G$ on $n$ vertices where $n$ is a power of $r$. This assumption is only made to simplify some expressions and, as we shall see, does not affect the asymptotic nature of our claims, since $i(F)$ and $i_k(F)$ are limits.

Recall that a {\em blowup} of $F$ is a graph obtained by replacing every vertex $v$ with an independent set $X_v$, and for any pair of vertices $u,v$ we have that if $uv \in E(F)$, then all possible $|X_u||X_v|$ edges between $X_u$ and $X_v$ are present, and if $uv \notin E(F)$, no such edge is present.
The {\em $t$-blowup} of $F$ is the blowup of $F$ where each blown-up part $X_v$ has order $t$.
A {\em fuzzy blowup} of $F$ is any graph obtained from a blowup of $F$ by possibly
adding some edges inside the blown-up parts.

We next define the {\em $s$-nested blowup} of $F$, denoted by $B_s(F)$. Defined inductively,
$B_1(F)=F$ and $B_s(F)$ is obtained by taking the lexicographic product of $F$ and $B_{s-1}(F)$.
That is, we take an $r^{s-1}$-blowup of $F$, and make each blown-up part induce a copy of $B_{s-1}(F)$.
Note that $B_s(F)$ has $r^s$ vertices.
It is easily observed that the number of induced copies of $F$ in $B_s(F)$ satisfies
$P(F,B_s(F)) \ge rP(F,B_{s-1}(F)) + (r^{s-1})^r$. This gives a geometric series for the lower bound from which we obtain
\begin{equation}\label{e:generic}
  P(F,B_s(F)) \ge r^{s-1} \left( \frac{r^{s(r-1)} - 1}{r^{r-1}-1}\right)\;.
\end{equation}
In turn, dividing by $\binom{r^s}{r}$ and taking $s$ to infinity yields the generic lower bound $i(F) \ge r!/(r^r-r)$ for inducibility, as observed by Pippenger and Golumbic \cite{PG-1975}.

We call $F$ a {\em fractilizer} if it holds for all $n=r^s$ that for all graphs $G$ on $n$ vertices,
$P(F,G)$ does not exceed the right hand side of \eqref{e:generic} and that the only graph which attains this bound is $B_s(F)$. In particular, if $F$ is a fractilizer, then $i(F) = r!/(r^r-r)$.

Let $B^t(F)$ denote the $t$-blowup of $F$ and observe that
$P(F,B^t(F)) \ge t^r$. Call $F$ {\em blowup maximal} if whenever a graph $G$ has $rt$ vertices and $P(F,G) \ge t^r$, then $G$ is a fuzzy $t$-blowup of $F$.

It has been proved by Fox, Huang, and Lee \cite{FHL-2017} that asymptotically almost surely,
$F \sim \G(r,\frac12)$ is a fractilizer. A similar result was obtained by the author in \cite{yuster-2019}.
In fact, in both papers, $F$ is asymptotically almost surely blowup maximal; see \cite{yuster-2019} Theorem 5.

We need one final additional property of $F$ which will guarantee that it is inducibility diverse.
Again, this property is also ``embedded'' in the constructions specified in the aforementioned papers, but as it is not explicitly stated there, we provide a proof here. Call $F$ {\em robust} if it is not a fuzzy blowup of any graph with at least two vertices and fewer than $r$ vertices. For example, $C_5$ is robust.

\begin{lemma}\label{l:robust}
	Asymptotically almost surely, $F \sim \G(r,\frac12)$ is a robust.
\end{lemma}
\begin{proof}
	Suppose that $V(F)=[r]$. Let $2 \le s \le r-1$. Suppose that $\P$ is a partition of $[r]$ into $s$ nonempty parts $V_1,\ldots,V_s$ where $|V_i| \ge |V_{i+1}|$ for $1 \le i \le s-1$. Let $X_{\P}$ denote the following event: For any two vertices $u,v$ in the same part $V_i$, and for any vertex $w \in V_j$ with $i \neq j$, either both $uw,vw$ are edges of $F$, or both are not. We will prove that the probability of $X_{\P}$ is small, so that the sum of these probabilities for each possible $\P$ is still $o_r(1)$, whence $F$ is robust with probability at least $1-o_r(1)$.
	
	Suppose that $|V_1|=q$ so we have that $2 \le q \le r-1$. For $X_{\cal P}$ to hold, the adjacencies of a particular vertex of $V_1$ with respect to the $r-q$ vertices not in $V_1$ determines the required adjacencies of any of the other $q-1$ vertices of $V_1$ with respect to the $r-q$ vertices not in $V_1$.
	In particular, $\Pr[X_\P] \le 2^{-(q-1)(r-q)}$.
	
	Consider first the case $q \ge 2\log_2 r$. Then there are fewer than $\binom{r}{q}$ choices for $V_1$
	and so at most $\binom{r}{q}(r-q)^{r-q}$ choices for such $\P$. Summing $\Pr[X_\P]$ for all such $\P$ gives for all sufficiently large $r$ at most
	\begin{equation}\label{e:robust-1}
	\sum_{q=2\log_2 r}^{r-1} \binom{r}{q}(r-q)^{r-q}  2^{-(q-1)(r-q)} < \frac{1}{r}\;.
	\end{equation}
	
	We may now assume that all parts of $\P$ have size at most $2\log_2r$. 
	Let $t$ be the number of singleton parts in $\P$. Consider first the case that $t \ge r-\sqrt{r}$. Let $u,v \in V_1$ (recall that $V_1$ is the largest part so is not a singleton). The probability that
	$u$ and $v$ have identical neighborhoods with respect to the singleton parts is $2^{-t}$.
	Notice that there are at most $\binom{r}{t}(r-t)^{r-t}$ choices for such $\P$.
	Summing $\Pr[X_\P]$ for all such $\P$ gives for all sufficiently large $r$ at most
	\begin{equation}\label{e:robust-2}
	\sum_{t=r-\sqrt{r}}^{r-2} \binom{r}{t}(r-t)^{r-t}  2^{-t} < \frac{1}{r}\;.
	\end{equation}
	
	We remain with the case that all parts of $\P$ have size at most $2\log_2r$ and there are at most
	$r-\sqrt{r}$ singleton parts. So, there are at least $\sqrt{r}$ vertices in non-singleton parts,
	hence there are $p \ge \sqrt{r}/2\log_2r$ non-singleton parts, say $V_1,\ldots,V_p$.
	Select two vertices $v_i,u_i$ from each $V_i$ for $1 \le i \le p$. The probability that $v_1,u_1$ see the same neighborhood in $V_2,\ldots,V_s$ is precisely $2^{|V_1|-r} \le 2^{-r/2}$.
	The probability that $v_i,u_i$ see the same neighborhood in $V_{i+1},\ldots,V_s$ (which is independent of whether $v_j,u_j$ see the same neighborhood in $V_{j+1},\ldots,V_s$) is precisely
	$2^{\sum_{j=1}^i|V_j|-r} \le 2^{2i\log_2 r-r} < 2^{-r/2}$ for all $1 \le i \le \sqrt{r}/2\log_2r \le p$.
	So,
	$$
	\Pr[X_\P] < (2^{-r/2})^{\sqrt{r}/2\log_2r} < 2^{-r^{4/3}}\;.
	$$
	As there are fewer than $r^r$ choices for $\P$, summing the last probability for all $\P$
	is at most $1/r$. Together with \eqref{e:robust-1} and \eqref{e:robust-2} we have that the probability that
	$F$ is robust is at least $1-3/r$ for all sufficiently large $r$.
\end{proof}

Since (trivially) almost all graphs have no isolated vertices, we have the following corollary.
\begin{corollary}\label{coro:properties}
	Asymptotically almost surely, $F \sim \G(r,\frac12)$ is a blowup-maximal robust fractilizer without isolated vertices.
\end{corollary}

\begin{proof}[Proof of Theorem \ref{t:fractil}]
	We shall prove that if $F$ is a blowup-maximal robust fractilizer without isolated vertices, then $F$ is inducibility diverse. Together with Corollary \ref{coro:properties}, this implies the theorem. So, suppose $F$ is such, and let $c=\omega(F)$ be the clique number of $F$.
	
	We prove two claims: (i) For each real $z < r!/(r^r-r)$ there exists $k \in \mathbb{N}$ such that $i_k(F) \ge z$. (ii) For each $k \in \mathbb{N}$ there exists a real $z < r!/(r^r-r)$ such that $i_k(F) \le z$. Clearly, the two claims together imply that $F$ is inducibility diverse. 
	
	For the first claim we do not really need the fact that $F$ is a fractilizer, yet it is easier to obtain explicit bounds if this is assumed. Suppose that $z < r!/(r^r-r)$.
	Recall from \eqref{e:generic} that
	$$
	\lim_{s \rightarrow \infty} \frac{P(F,B_s(F))}{\binom{r^s}{r}} = \frac{r!}{r^r-r}\;.
	$$
	So, let $s$ be the smallest integer such that $p(F,B_s(F))=P(F,B_s(F))/\binom{r^s}{r} \ge z$.
	But further observe that since the nested blowup is defined via lexicographic product,
	we have that $\omega(B_s(F))=c^s$, so let $k=c^s$. Now consider $n$-th blowup of $B_s(F)$,
	namely $G=B^n(B_s(F))$. Notice that $\omega(G)=k$ and that $p(F,G)=p(F,B_s(F))(1-o_n(1)) \ge z(1-o_n(1))$. Thus, we have that $i_k(F) \ge z$, as required.
	
	We prove the second claim by induction on $k$ and notice that we can start the base case at any $k$, so we start at $k=c$. We claim that $i_{c+1}(F)=r!/r^r < r!/(r^r-r)$. Indeed, observe that a $t$-blowup of $F$ has
	$t^r$ induced copies of $F$, so dividing by $\binom{tr}{r}$ already gives $i_{c+1}(F) \ge r!/r^r$.
	On the other hand, since $F$ is blowup maximal and robust, any graph $G$ with $rt$ vertices and more than $t^r$ induced copies of $F$ must be a fuzzy $t$-blowup of $F$ and, moreover, some blown-up part
	contains an induced copy of $F$. Suppose this is the blown-up part $X_v$ of vertex $v$ and let $K$ be a $c$-clique of an induced copy of $F$ in $X_v$. Let $u$ be a neighbor of $v$ in $F$ (recall that $F$ has no isolated vertices). Take any vertex $y$ of $X_u$ and notice that $y$ is connected in $G$ to all vertices of $K$, so $G$ contains a $(c+1)$-clique. It thus immediately follows that $i_{c+1}(F) \le r!/r^r$.
	
	Next, suppose the claim holds for $k \ge c$ and we shall prove that it holds for $k+1$ as well.
	Let $z^* > i_k(F)$ where $z^* < r!/(r^r-r)$ exists by the induction hypothesis.
	Let $z$ be chosen such that
	$$
	\frac{r!+rz^*}{r^r} < z < \frac{r!}{r^r-r}\;.
	$$
	Notice that $z$ exists since
	$$
	\frac{r!+rz^*}{r^r} < \frac{r!+r(r!/(r^r-r))}{r^r} = \frac{r!}{r^r-r}\;.
	$$
	We shall prove that $i_{k+1}(F) \le z$.
	Indeed, take any graph $G$ with $r^s$ vertices
	for which $p(F,G) \ge z$. Let $t=r^{s-1}$ and observe that for all $s$ sufficiently large,
	\begin{equation}\label{e:pfg}
	P(F,G) = \binom{rt}{r}p(F,G) \ge \binom{rt}{r}z > \binom{rt}{r}\frac{r!+rz^*}{r^r} = (1-o_s(1))\frac{r^rt^r}{r!} \cdot \frac{r!+rz^*}{r^r} > t^r\;.
	\end{equation}
	But since $F$ is blowup maximal, this implies that $G$ is a fuzzy $t$-blowup of $F$.
	So, let the parts of this fuzzy $t$-blowup be $\{X_v\,|\, v \in V(F)\}$ observing that
	$|X_v|=r^{s-1}$ for all $v \in V(F)$. Since $F$ is robust, it must be that each induced copy of
	$F$ in $G$ is either a transversal having precisely one vertex in each $X_v$, or else it is entirely
	contained in some $X_v$. So, the number of induced copies of $F$ in $G$ that are entirely contained in some $X_v$ is $P(F,G)-t^r$. By pigeonhole, there is some $X_v$ which contains at least
	$(P(F,G)-t^r)/r$ induced copies of $F$, so we pick one. Let $G'=G[X_v]$ be the induced subgraph of $G$ on vertex set $X_v$. We have by \eqref{e:pfg} 
	$$
	P(F,G') \ge \frac{P(F,G)-t^r}{r} \ge (1-o_s(1)) \frac{\frac{t^r}{r!}(r!+rz^*)-t^r}{r} = (1-o_s(1)) \binom{t}{r}z^*\;.
	$$
	Hence, $p(F,G') \ge (1-o_s(1))z^*$. But since $z^* > i_k(F)$ we have that for $s$ sufficiently large,
	$G'$ must contain a clique $K$ or order $k$. Let $u$ be a neighbor of $v$ in $F$. Take any vertex $y$ of $X_u$ and notice that $y$ is connected in $G$ to all vertices of $K$, so $G$ contains a $(k+1)$-clique. It follows that $i_{k+1}(F) \le z$. 
\end{proof}

\section{Inducibility of symmetrizable families}\label{sec:symmetrizable}

In this section we prove Theorem \ref{t:symmetrizable}. Throughout this section, $F$ denotes a symmetrizable family of graphs on $r$ vertices and $H$ is a graph with $\chi(H)=k \ge 3$.
We shall prove that $i_H(F)=i_k(F)$ and that $i_k(F)$ (whence $i_H(F)$) is a attained by a complete $\ell$-partite graphon where $\ell < k$.

Before we prove the theorem, let us first consider the ``toy-case'' where $H=K_k$, and use Zykov symmetrization to show that among all $K_k$-free graphs $G$ with $n$ vertices, $P(F,G)$ is maximized when $G$ is a complete partite graph.
For a vertex $v$ of $G$, let $s(v)$ denote the number of induced $F$-copies containing $v$ (recall that $F$ is a family, so by induced $F$-copy, we mean and induced subgraph of $G$ that is isomorphic to an element of $F$). If $u,v$ are nonadjacent and $s(u) < s(v)$, then replace
$u$ with a copy of $v$ (keeping $v$ and its clone nonadjacent).
As $F$ is symmetrizable, $P(F,G)$ increases by $s(v)-s(u)$ and clearly the new graph is still $K_k$-free. To see the former, just notice that if an $F$-copy containing both $u$ and $v$ existed before symmetrization, then it is also an $F$-copy after symmetrization (albeit that copy may now correspond to another element of $F$).
Proceeding in this manner, once done, we may assume that any two nonadjacent vertices have the same $s(.)$ value. Two vertices $u$ and $v$ are {\em twins} if they have identical neighborhoods (in particular, twins are nonadjacent), and notice that twinhood is an equivalence relation. Let $X$ and $Y$ be two distinct equivalence classes and observe that we must have that either each vertex of $X$ is adjacent to each vertex of $Y$ or that $X \cup Y$ is an independent set. Assume the latter. We can ``absorb'' $Y$ into $X$ as follows. Let $y \in Y$ and replace $y$ with a copy of (any) $v \in X$. As $F$ is symmetrizable, the number of induced copies of $F$ in $G$ does not change, yet the graph remains $K_k$-free,
$|X|$ increases by $1$ and $|Y|$ decreases by $1$. Proceeding in this manner, once done, we may assume that any for two equivalence classes $X$ and $Y$, each vertex of $X$ is adjacent to each vertex of $Y$; in other words, $G$ is complete partite.

To prove the theorem for general $H$, we shall need Szemer\'edi's Regularity Lemma \cite{szemeredi-1978} which we briefly recap.
We say that a pair of disjoint nonempty vertex sets $A$, $B$ of a graph are {\em $\varepsilon$-regular}
if for all $X \subseteq A$ and $Y \subseteq B$ with $|X|\ge \varepsilon|A|$ and $|Y|\ge \varepsilon |B|$,
$$
\left|\frac{e(X, Y)}{|X||Y|}-\frac{e(A, B)}{|A||B|}\right| \leq \varepsilon
$$
where $e(U,V)$ denotes the number of edges with an endpoint in $U$ and an endpoint in $V$.
An $\varepsilon$-regular partition of a graph $G$ is a partition of $V(G)$ into sets $V_1, \ldots , V_\ell$ such that $\ell \ge \varepsilon^{-1}$, $||V_i|-|V_j||\le 1$ for all $i,j\in[\ell]$,
and all but $\varepsilon\ell^2$ pairs $V_i,V_j$ are $\varepsilon$-regular.
Szemer\'edi's regularity lemma states that for every $\varepsilon > 0$ there exists $K(\varepsilon)$ such that every graph with at least $\varepsilon^{-1}$ vertices has an $\varepsilon$-regular partition with at most $K(\varepsilon)$ parts.
One of the canonical applications of the regularity lemma is its use to prove a corresponding {\em removal lemma} \cite{EFR-1986}: Suppose that $\delta > 0$ is some small given constant and that
$G$ is a large $n$-vertex $H$-free graph with $\chi(H)=k \ge 3$. Let $\varepsilon=\varepsilon(\delta,H)$ be sufficiently small, and consider an $\varepsilon$-regular partition of $G$ into parts
$V_1,\ldots,V_\ell$ with the properties guaranteed by the regularity lemma.
Modify $G$ by removing (i) edges with both endpoints in the same part (ii) edges between two non $\varepsilon$-regular pairs (iii) edges between pairs for which $e(V_i,V_j)/|V_i||V_j| \le \delta$.
Observe that the number of removed edges is smaller than
$$
\ell (n^2/\ell^2) + \varepsilon \ell^2 (n^2/\ell^2) + \ell^2(n^2/\ell^2)\delta \le n^2(\delta+2\varepsilon) < 2\delta n^2\;.
$$
It is a standard embedding argument to prove that the new graph $G'$ has no $K_k$. Indeed, if it had one,
then its vertices must span $k$ distinct vertex parts, and all $\binom{k}{2}$ of the pairs of participating
parts are $\varepsilon$-regular, and have edge density at least $\delta$. But now we can use the
$\varepsilon$-regularity and the assumed density of these pairs to embed in $G'$ a compete $k$-partite graph with $|V(H)|$ vertices in each part, and in particular, embed a copy of $H$ in $G'$ which is impossible as $G'$ is a subgraph of $G$.
One must be careful though, as in $G'$ there may be fewer $F$-copies as there were in $G$.
However, observe that each pair of vertices of $G$ trivially participates in fewer than $n^{r-2}$ $F$-copies,
thus $|P(F,G)-P(F,G')| \le 2\delta n^r$.
We therefore obtain the following corollary.
\begin{corollary}\label{coro:reg}
	Let $\gamma > 0$. For all $n$-sufficiently large, if $G$ is an $H$-free graph with $n$ vertices, then $G$ has a spanning subgraph $G'$ that is $K_k$-free and such that $|P(F,G)-P(F,G')| \le \gamma n^r$.
\end{corollary}

\begin{proof}[Proof of Theorem \ref{t:symmetrizable}]
	We first prove that $i_k(F) \le i_H(F)$.
	Fix $n$, and consider a graph $G$ that attains $i_k(F,n)$. By the discussion above, we can assume that $G$ is complete partite, i.e. complete $\ell$-partite for some $\ell < k$. As $\chi(G)=\ell$,
	we have that $G$ is also $H$-free. Hence, $i_H(F,n) \ge i_k(F,n)$. In particular, at the limit, we have
	$i_k(F) \le i_H(F)$.
	
	We next prove that $i_k(F) \ge i_H(F)$.  Let $\gamma > 0$ and let $n$ be sufficiently large such that
	Corollary \ref{coro:properties} applies. Consider a graph $G$ that attains $i_H(F,n)$.
	Then, by Corollary \ref{coro:reg}, $G$ has a spanning subgraph $G'$ that is $K_k$-free and
	$P(F,G') \ge P(F,G) - \gamma n^r$.  Hence,
	$$
	i_k(F,n) \ge p(F,G')  \ge i_H(F,n) - \frac{\gamma n^r}{\binom{n}{r}} \ge i_H(F,n)-\gamma r^r\;.
	$$
	In particular, at the limit, we have that $i_k(F) \ge i_H(F)- \gamma r^r$.
	As this holds true for every $\gamma > 0$, we have that $i_k(F) \ge i_H(F)$.
	
	Since the graphon attaining $i_k(F)$ is $K_k$-free, and since by the aforementioned Zykov symmetrization detailed in the beginning of this section, this graphon is a limit of complete partite graphs, it must be that the graphon attaining $i_H(F)=i_k(F)$ is complete $\ell$-partite for some $\ell < k$.
\end{proof}

\section{Tur\'an graphs}\label{sec:turan}
\subsection{Polynomial optimization}

For a complete partite graph $F=K_{a_1,\ldots,a_r}$ and for $k \ge 3$,
we know by Theorem \ref{t:symmetrizable} that $i_k(F)$ is attained by some complete $\ell$-partite graphon
$W_{F,k}$ where $\ell < k$. Since for a given graphon $W[x_1,\ldots,x_{k-1}]$, the induced density of $F$ in $W$
is determined by $F$ and $x_1,\ldots,x_{k-1}$, and since this density is a symmetric polynomial $P_{F,k} = P(x_1,\ldots,x_{k-1})$ of degree $|V(F)|$, we have that determining $i_k(F)$ amounts to solving a polynomial optimization problem: maximizing $P_{F,k}$ over $\Delta^{k-2}$ (in particular, $\ell$ is the number of nonzero coordinates of such a solution). While generally challenging,
for some $F$ we are able to solve this optimization problem, and for some others, determine
certain distinguishable characteristics of a solution. In the remaining of this section we consider Tur\'an graphs and in the following section we consider other complete partite graphs.

\begin{example}
Suppose that $F=K_{2,1,1}$ and $k=5$.
Then the graphon has (at most) four parts of sizes $x_1,x_2,x_3,x_4$ where $(x_1,x_2,x_3,x_4) \in \Delta^3$.
The probability that a random placement of four points in $[0,1]$ will span a $K_{2,1,1}$ in this graphon is therefore
\begin{align*}
P_{F,5} & =  12x_1^2(x_2x_3+x_2x_4+x_3x_4) + 12x_2^2(x_1x_3+x_1x_4+x_3x_4)\\
&  + 12x_3^2(x_1x_2+x_1x_4+x_2x_4) + 12x_4^2(x_1x_2+x_1x_3+x_2x_3)\;.
\end{align*}
It is not difficult to show that $P_{F,5}$ is maximized at $(\frac{1}{4},\frac{1}{4},\frac{1}{4},\frac{1}{4})$ where we have $i_5(K_{2,1,1})=\frac{9}{16}$ and hence by Theorem \ref{t:symmetrizable}, $i_H(K_{2,1,1})=\frac{9}{16}$ for every graph $H$ with $\chi(H)=5$.
\end{example}

\subsection{Highly dense Tur\'an graphs and complete equipartite graphs}
In this subsection we prove Theorem \ref{t:turan}, that is, we determine the (unique) solution
maximizing $P_{F,k}$ for $F=T(s,r)$ where $2 \le r < s \le 3r+1$. 
Clearly, we may assume that $k > r$ as otherwise we trivially have $i_k(T(s,r))=0$.
Furthermore, if $k > r$, then $i_k(T(s,r)) > 0$.

\begin{proof}[Proof of Theorem \ref{t:turan}]
	Consider the Tur\'an graph $F=T(s,r)$ where $2 \le r < s$ (at this point we do not assume that $s \le 3r+1$, we shall require this later). For notation clarity, let $s=pr+q$ where $0 \le q < r$.
	That it, $F$ has $q$ parts of size $p+1$ and $r-q$ parts of size $p$.
	
	Let ${\bf x} = (x_1,\ldots,x_{k-1}) \in \arg\max_{{\bf x} \in \Delta^{k-2}} P_{F,k}$ and consider the corresponding graphon $W[{\bf x}]$. Let $\ell \le k-1$ be the number of nonzero entries in ${\bf x}$
	(so $W$ is, in fact, a complete $\ell$-partite graphon). Also assume by symmetry of $P_{F,k}$
	that $x_i \ge x_{i+1}$ for $i=1,\ldots,k-2$. We shall first prove that $W$ is equipartite.
	
	\begin{lemma}\label{l:equal}
		$x_i=1/\ell$ for $i=1,\ldots,\ell$.
	\end{lemma}
	\begin{proof}
		Assume that the statement is false, and let $i < \ell$ be the smallest index such that
		$x_i > x_{i+1} > 0$. Let $a$ be a variable, replace $x_i$ with $a(x_i+x_{i+1})$
		and replace $x_{i+1}$ with $(1-a)(x_i+x_{i+1})$. Observe that for every $a \in [0,1]$,
		the point $p(a)=(x_1,\ldots,x_{i-1},a(x_i+x_{i+1}),(1-a)(x_i+x_{i+1}),x_{i+2},\ldots,x_{k-1})$ is still on
		$\Delta^{k-2}$ and that for $a=x_i/(x_i+x_{i+1}) \notin \{0,\frac{1}{2},1\}$ the corresponding point is $p(a)={\bf x} \in 
		\arg\max_{{\bf x} \in \Delta^{k-2}} P_{F,k}$. We shall derive a contradiction by showing that there
		is another point $a \in [0,1]$ for which $P_{F,k}(p(a)) > P_{F,k}({\bf x})$.
		
		Fixing $x_1,\ldots,x_{k-1}$ and letting $a$ vary, consider the univariate polynomial $Q(a)=P_{F,k}(p(a))$
		in the entire real line. (For example, if $F=T(5,3)$ and $k=4$ we have that
		$P_{F,4}=30x_1x_2x_3(x_1x_2+x_1x_3+x_2x_3)$ and if $x_1 > x_2 > 0$ we have
		$Q(a) = 30((1-a)^2 a^2 x_3 (x_1 + x_2)^4 + a(1-a)(x_1 + x_2)^3x_3^2)$.)
		Some general observable facts on $Q(a)$ follow.
		
		First notice that $Q(a)$ is symmetric around
		$\frac{1}{2}$ since $Q(\frac{1}{2}+z)=Q(\frac{1}{2}-z)$ for every $z \in {\mathbb R}$.
		Hence, $Q'(\frac{1}{2})=0$.
		
		A copy of $F$ in the corresponding graphon has six possibilities: (i) it has $p+1$ vertices in part $i$ and $p+1$ vertices in part $i+1$, (ii) it has $p+1$ vertices in part $i$ and $p$ vertices in part $i+1$ or vice versa,
		(iii) it has $p+1$ vertices in part $i$ and no vertex in part $i+1$ or vice versa,
		(iv) it has $p$ vertices in 
		part $i$ and $p$ vertices in part $i+1$, (v) it has $p$ vertices in part $i$ and no vertex in part $i+1$ or vice versa, (vi) it has no vertex in part $i$ and no vertex in part $i+1$.
		
		In what follows $C_1,\ldots,C_6$ are polynomials depending each only on $x_1,\ldots,x_{k-1}$
		and each is a sum of positive monomials, so each is positive at $(x_1,\ldots,x_{k-1})$.
		The contribution of copies of case (vi) to $Q(a)$ is some constant $C_6$.
		The contribution of copies of case (v) to $Q(a)$ is of the form
		$a^p(x_i+x_{i+1})+(1-a)^p(x_i+x_{i+1})$ multiplied by a constant so of the form
		$(a^p+(1-a)^p)C_5$.
		The contribution of copies of case (iv) to $Q(a)$ is of the form
		$a^p(x_i+x_{i+1})(1-a)^p(x_i+x_{i+1})$ multiplied by a constant, so is of the form $a^p(1-a)^p C_4$.
		The contribution of copies of case (iii) to $Q(a)$ is of the form
		$a^{p+1}(x_i+x_{i+1})+(1-a)^{p+1}(x_i+x_{i+1})$ multiplied by a constant, so is of the form $(a^{p+1}+(1-a)^{p+1})C_3$.
		The contribution of copies of case (ii) to $Q(a)$ is of the form
		$a^{p+1}(x_i+x_{i+1})(1-a)^p(x_i+x_{i+1})+a^p(x_i+x_{i+1})(1-a)^{p+1}(x_i+x_{i+1})$ multiplied by a constant, so is of the form $a^p(1-a)^p C_2$.
		The contribution of copies of case (i) to $Q(a)$ is of the form
		$a^{p+1}(x_i+x_{i+1})(1-a)^{p+1}(x_i+x_{i+1})$ multiplied by a constant, so is of the form
		$a^{p+1}(1-a)^{p+1}C_1$.	
		
		Consider first the case $p=1$. In this case, $Q(a)$ has degree $4$ and, moreover, the coefficient of $a^4$ is $C_1 > 0$.
		Therefore, $\lim_{a\rightarrow \pm \infty}P(a) = +\infty$.
		Hence, $Q(a)$ has either one extremal point at $a=\frac{1}{2}$ or three extremal points:
		a local maximum at $a=\frac{1}{2}$ and two minima. In both cases, we have that
		$\arg \max_{a \in [0,1]}Q(a) \in \{0,\frac{1}{2},1\}$, a contradiction.
		
		Consider next the case $p=2$. In this case, $Q(a)$ has degree $6$ and, moreover, the coefficient of $a^6$ is $-C_1 < 0$. Therefore, $\lim_{a\rightarrow \pm \infty}P(a) = -\infty$.
		Now,
		$$
		Q'(a) = (4a-2)C_5+(4a^3-6a^2+2a)(C_4+C_2)+(6a-3)C_3+(-6a^4+15a^4-12a^3+3a^2)C_1
		$$
		so $Q'(0) = -2C_5-3C_3 < 0$ and similarly $Q'(1)=-Q(0) > 0$.
		So $Q$ is increasing at $a=1$ and thus has a local maximum at a point larger than $1$ and similarly,
		it has a local maximum at a point smaller than $0$. So either $Q$ has five local extrema, which means
		that $\frac{1}{2}$ is the remaining local maximum or else it has three local extrema, which means that $\frac{1}{2}$ is the remaining local extremum (minimum, in fact). 
		In both cases, we have that
		$\arg \max_{a \in [0,1]}Q(a) \in \{0,\frac{1}{2},1\}$, a contradiction.
		
		Consider next the case $p=3$ and $q=0$ (so all parts of $F$ have size $3$).
		Here we have that the only possible cases are (iv), (v) and (vi) which means that
		$$
		Q(a)=C_6+(3a^2-3a+1)C_5+(a^3 - 3 a^4 + 3 a^5 - a^6)C_4\;.
		$$
		Here again we have $\lim_{a\rightarrow \pm \infty}P(a) = -\infty$ and
		$$
		Q'(a) = (6a-3)C_5+(3a^2-12a^3+15a^4-6a^5)C_4
		$$
		so $Q'(0) = -3C_5 < 0$ and similarly $Q'(1)=-Q(0) > 0$.
		Thus, we arrive at precisely the same situation as in the previous paragraph for the case $p=2$, which, recall, is a contradiction.
		
		Consider finally the case $p=3$ and $q=1$ (so all parts of $F$ but one have size $3$ and one part has size $4$).
		Here we have that all cases are (ii), (iii), (iv), (v) and (vi) are possible, but Case (i) is impossible. This means that
		$$
		Q(a)=C_6+(3a^2-3a+1)C_5+(a^3 - 3 a^4 + 3 a^5 - a^6)(C_4+C_2)+(2a^4 - 4a^3+6a^2-4a+1)C_3\;.
		$$
		Here again we have $\lim_{a\rightarrow \pm \infty}P(a) = -\infty$ and $Q'(0) = -3C_5-4C_3 < 0$ and similarly $Q'(1)=-Q(0) > 0$.
		Thus, we arrive at precisely the same situation as in the case $p=2$, which, recall, is a contradiction.
		
		Since the statement holds for $p=1,2$ and for $(p,q)=(3,0)$ and $(p,q)=(3,1)$, this covers all possibilities as we assume $s \le 3r+1$.
	\end{proof}

Having proved Lemma \ref{l:equal}, it remains to optimize $\ell$. Let $g(\ell)$ denote the induced density
of $T(s,r)$ in $W[\frac{1}{\ell},\ldots,\frac{1}{\ell}]$. It is an easy exercise to see that
\begin{equation}\label{e:g}
  g(\ell) = \frac{\ell!(pr+q)!}{(\ell-r)!(r-q)!q!(p!)^r(p+1)^q\ell^{pr+q}}\;.
\end{equation}

\begin{lemma}\label{l:unimodal}
	The discrete function $g(\ell)$ for $\ell \ge r$ is unimodal; staring from $\ell=r$, $g(\ell)$ strictly increases until it attains a maximum at some unique point $t$ (possibly $t=r$), and then strictly decreases.
\end{lemma}
\begin{proof}
Consider the ratio
$g(\ell)/g(\ell-1)$ which is
$$
f(\ell) = \frac{g(\ell)}{g(\ell-1)} = \frac{\ell}{\ell-r}\left(1-\frac{1}{\ell}\right)^{pr+q}\;.
$$
Viewing $f(x)$ as a real function in $(r,\infty)$ and considering $f'$, we see that
$\lim_{x \rightarrow r^+} = \infty$, $f$ has a horizontal asymptote at $1$, $f(x) < 1$ for all
sufficiently large $x$ (since $pr+q = s \ge r+1$) and $f'(x)=0$ only at $x=r(pr+q-1)/(q+(p-1)r)$. Hence there is a single
point $x_0$ at which $f(x_0)=1$, $f(x)  > 1$ for $x < x_0$ and $f(x) < 1$ for $x > x_0$.
Hence, the required $t$ is just the largest $\ell$  for which $f(\ell) > 1$ (note that $t=r$
if $f(r+1) < 1$).
\end{proof}

Having proved Lemmas \ref{l:equal} and \ref{l:unimodal} we have that if $k \le t+1$, then
$i_k(T(s,r))$ is attained by $W[\frac{1}{k-1},\ldots,\frac{1}{k-1}]$ and furthermore $i_k(T(s,r))=g(k-1)$,
and if $k > t+1$, then $i_k(T(s,r))$ is attained by $W[\frac{1}{t},\ldots,\frac{1}{t}]$ and furthermore $i_k(T(s,r))=g(t)$.
In particular, $i(T(s,r))$ is attained by $W[\frac{1}{t},\ldots,\frac{1}{t}]$ and furthermore $i(T(s,r))=g(t)$.
This completes the proof of Theorem \ref{t:turan}.
\end{proof}

\subsection{All Tur\'an graphs on at most $14$ vertices}
In this subsection we prove Theorem \ref{t:le14}.

\begin{table}[p!t]
	\tiny
	\centering
	\begin{tabular}{c||c|c|c|c}
		\small{graph} & \small{$t$} & \small{numerator} & \small{denominator} & \small{reference}\\
		\hline
		$T(3,2)$ & $2$ & $3$ & $4$ & \cite{goodman-1959}\\
		\hline
		$T(4,2)$ & $2$ & $3$ & $8$ & \cite{BS-1994}\\
		\hline
		$T(4,3)$ & $5$ & $72$ & $125$ & \cite{hirst-2014}\\
		\hline
		$T(5,2)$ & $2$ & $5$ & $8$ & \cite{BS-1994}\\
		\hline
		$T(5,3)$ & $3$ & $10$ & $27$ & \cite{BS-1994} \cite{PST-2019}\\
		\hline
		$T(5,4)$ & $8$ & $525$ & $1024$ & \cite{LPSS-2023} \cite{LMR-2023}\\
		\hline
		$T(6,2)$ & $2$ & $5$ & $16$ & \cite{BS-1994}\\
		\hline
		$T(6,3)$ & $3$ & $10$ & $81$ & \cite{BS-1994}\\
		\hline
		$T(6,4)$ & $6$ & $25$ & $72$ & here\\
		\hline
		$T(6,5)$ & $13$ & $178200$ & $371293$ & \cite{LMR-2023}\\
		\hline
		$T(7,2)$ & $2$ & $35$ & $64$ & \cite{BS-1994}\\
		\hline
		$T(7,3)$ & $3$ & $70$ & $243$ & here\\
		\hline
		$T(7,4)$ & $5$ & $504$ & $3125$ & here\\
		\hline
		$T(7,5)$ & $8$ & $11025$ & $32768$ & here\\
		\hline
		$T(7,6)$ & $19$ & $21591360$ & $47045881$ & \cite{LMR-2023}\\
		\hline
		$T(8,2)$ & $2$ & $35$ & $128$ & \cite{BS-1994}\\
		\hline
		$T(8,3)$ & $3$ & $560$ & $2187$ & here\\
		\hline
		$T(8,4)$ & $4$ & $315$ & $8192$ & \cite{BS-1994}\\
		\hline
		$T(8,5)$ & $7$ & $21600$ & $117649$ & here\\
		\hline
		$T(8,6)$ & $11$ & $6350400$ & $19487171$ & here\\
		\hline
		$T(8,7)$ & $25$ & $542691072$ & $1220703125$ & \cite{LMR-2023}\\
		\hline
		$T(9,2)$ & $2$ & $63$ & $128$ & \cite{BS-1994}\\
		\hline
		$T(9,3)$ & $3$ & $560$ & $6561$ & \cite{BS-1994}\\
		\hline
		$T(9,4)$ & $4$ & $945$ & $8192$ & here\\
		\hline
		$T(9,5)$ & $6$ & $175$ & $2592$ & here\\
		\hline
		$T(9,6)$ & $9$ & $313600$ & $1594323$ & here\\
		\hline
		$T(9,7)$ & $15$ & $224224$ & $703125$ & here\\
		\hline
		$T(9,8)$ & $33$ & $837724160$ & $1929229929$ & \cite{LMR-2023}\\
		\hline
		$T(10,2)$ & $2$ & $63$ & $256$ & \cite{BS-1994}\\
		\hline
		$T(10,3)$ & $3$ & $1400$ & $6561$ & \cite{BS-1994}\\
		\hline
		$T(10,4)$ & $4$ & $4725$ & $32768$ & here\\
		\hline
		$T(10,5)$ & $5$ & $4536$ & $390625$ & \cite{BS-1994}\\
		\hline
		$T(10,6)$ & $8$ & $1488375$ & $16777216$ & here\\
		\hline
		$T(10,7)$ & $12$ & $67375$ & $331776$ & here\\
		\hline
		$T(10,8)$ & $19$ & $101047564800$ & $322687697779$ & here\\
		\hline
		$T(10,9)$ & $42$ & $796235375$ & $1867795524$ & \cite{LMR-2023}\\
		\hline
		$T(11,2)$ & $2$ & $231$ & $512$ & \cite{BS-1994}\\
		\hline
		$T(11,3)$ & $3$ & $3850$ & $19683$ & \cite{BS-1994}\\
		\hline
		$T(11,4)$ & $4$ & $5775$ & $65536$ & here\\
		\hline
		$T(11,5)$ & $5$ & $16632$ & $390625$ & here\\
		\hline
		$T(11,6)$ & $7$ & $1069200$ & $40353607$ & here\\
		\hline
		$T(11,7)$ & $10$ & $130977$ & $1250000$ & here\\
		\hline
		$T(11,8)$ & $15$ & $19731712$ & $94921875$ & here\\
		\hline
		$T(11,9)$ & $24$ & $157309075$ & $509607936$ & here\\
		\hline
		$T(11,10)$ & $51$ & $617338863680000$ & $1469659553427321$ & \cite{LMR-2023}\\
		\hline
		$T(12,2)$ & $2$ & $231$ & $1024$ & \cite{BS-1994}\\
		\hline
		$T(12,3)$ & $3$ & $3850$ & $59049$ & \cite{BS-1994}\\
		\hline
		$T(12,4)$ & $4$ & $5775$ & $262144$ & \cite{BS-1994}\\
		\hline
		$T(12,5)$ & $5$ & $133056$ & $1953125$ & here\\
		\hline
		$T(12,6)$ & $7$ & $1069200$ & $282475249$ & here\\
		\hline
		$T(12,7)$ & $9$ & $1724800$ & $43046721$ & here\\
		\hline
		$T(12,8)$ & $12$ & $3705625$ & $31850496$ & here\\
		\hline
		$T(12,9)$ & $18$ & $163788625$ & $774840978$ & here\\
		\hline
		$T(12,10)$ & $29$ & $3721968169920000$ & $12200509765705829$ & here\\
		\hline
		$T(12,11)$ & $62$ & $21091312901233575$ & $50816953792809662$ & \cite{LMR-2023}\\
		\hline
		$T(13,2)$ & $2$ & $429$ & $1024$ & \cite{BS-1994}\\
		\hline
		$T(13,3)$ & $3$ & $10010$ & $59049$ & \cite{BS-1994}\\
		\hline
		$T(13,4)$ & $4$ & $75075$ & $1048576$ & \cite{BS-1994}\\
		\hline
		$T(13,5)$ & $5$ & $576576$ & $9765625$ & here\\
		\hline
		$T(13,6)$ & $6$ & $25025$ & $1679616$ & here\\
		\hline
		$T(13,7)$ & $8$ & $42567525$ & $4294967296$ & here\\
		\hline
		$T(13,8)$ & $11$ & $14859936000$ & $285311670611$ & here\\
		\hline
		$T(13,9)$ & $15$ & $897792896$ & $7119140625$ & here\\
		\hline
		$T(13,10)$ & $22$ & $243737778375$ & $1141246682444$ & here\\
		\hline
		$T(13,11)$ & $35$ & $233296472927232$ & $772393258984375$ & here\\
		\hline
		$T(13,12)$ & $74$ & $2704936173225986700$ & $6582952005840035281$ & \cite{LMR-2023}\\
		\hline
		$T(14,2)$ & $2$ & $429$ & $2048$ & \cite{BS-1994}\\
		\hline
		$T(14,3)$ & $3$ & $28028$ & $177147$ & \cite{BS-1994}\\
		\hline
		$T(14,4)$ & $4$ & $1576575$ & $16777216$ & \cite{BS-1994}\\
		\hline
		$T(14,5)$ & $5$ & $1345344$ & $48828125$ & here\\
		\hline
		$T(14,6)$ & $6$ & $875875$ & $30233088$ & here\\
		\hline
		$T(14,7)$ & $8$ & $42567525$ & $34359738368$ & here\\
		\hline
		$T(14,8)$ & $10$ & $107270163$ & $6250000000$ & here\\
		\hline
		$T(14,9)$ & $13$ & $1452272976000$ & $23298085122481$ & here\\
		\hline
		$T(14,10)$ & $18$ & $14904764875$ & $111577100832$ & here\\
		\hline
		$T(14,11)$ & $26$ & $20068705243125$ & $93192340489924$ & here\\
		\hline
		$T(14,12)$ & $41$ & $277115986307736576000$ & $925103102315013629321$ & here\\
		\hline
		$T(14,13)$ & $86$ & $2800242205096869658125$ & $6873056497129163140972$ & \cite{LMR-2023}	
	\end{tabular}
	\caption{The inducibility values and the attaining equipartite graphon on $t$ parts for Tur\'an graphs on at most $14$ vertices.}
	\label{table:le14} 
\end{table} 

For each $3 \le s \le 14$ and each $2 \le r < s$ (the case $r=s$ is the case of complete graphs which is
given in Corollary \ref{coro:erdos}), Table \ref{table:le14} list $i(T(s,r))$ in the form of a numerator and denominator. In all cases, $i(T(s,r))$ is attained uniquely by $W[\frac{1}{t},\ldots,\frac{1}{t}]$ where $t$ is given in the table. In all cases, we list the reference of the first result which proves that value. To obtain $i_k(T(s,r))$ from that table, we may proceed as follows.
If $k \le r$, then trivially, $i_k(T(s,r))=0$. If $r < k \le t$, then $i_k(T(s,r))$ is uniquely attained by
$W[\frac{1}{k-1},\ldots,\frac{1}{k-1}]$, and if $k > t$, then $i_k(T(s,r))$ is uniquely attained by
$W[\frac{1}{t},\ldots,\frac{1}{t}]$. To obtain the actual value of $i_k(T(s,r))$ in the case that $r < k \le t$,
we note that in all such listed cases, it always holds that the conditions of Theorem \ref{t:turan} are met and so $i_k(T(s,r))=g(k-1)$ using \eqref{e:g}.

To verify the correctness of values and the table, we only need to consider the cases in which the conditions in Theorem \ref{t:turan} do not apply, as when they do, they are a consequence of Theorem \ref{t:turan}. For the bipartite cases when $r=2$, the result of Brown and Sidorenko \cite{BS-1994} gives
that $W[\frac{1}{2},\frac{1}{2}]$ is the unique graphon for $T(s,2)$. In particular,
$i_k(T(s,2)$ is also uniquely attained by $W[\frac{1}{2},\frac{1}{2}]$ for all $k \ge 3$.
The only other graphs in this table which the conditions in Theorem \ref{t:turan} do not apply
are $T(11,3)$, $T(12,3)$, $T(13,3)$, $T(14,3)$, $T(14,4)$. However, in all of these cases,
the condition in \eqref{e:BS} applies and we have that $W[\frac{1}{r},\ldots,\frac{1}{r}]$ is the unique graphon for $T(s,r)$ in these cases. In particular,
$i_k(T(s,r)$ is also uniquely attained by $W[\frac{1}{r},\ldots,\frac{1}{r}]$ for all $k \ge r+1$.
The smallest Tur\'an graph which is not bipartite and not covered by Theorem \ref{t:turan} nor
\eqref{e:BS} is $T(15,4)$.  \qed

\subsection{Proof of Proposition \ref{prop:turan}}

We need to prove that $i_{r+1}(T(s,r))$ is attained by $W[\frac{1}{r},\ldots,\frac{1}{r}]$ for all $2 \le r \le s$. We shall use the same arguments and notation as in the proof of Lemma 
\ref{l:equal} (where here $\ell=r$). Recall that $s=pr+q$ where $0 \le q < r$
so that $T(s,r)$ has $q$ parts of size $p+1$ and $r-q$ parts of size $p$.

As any copy of $F$ in the graphon $W[x_1,\ldots,x_r]$ must contain vertices in each part,
when we consider $Q(a)$ in Lemma \ref{l:equal} only the cases (i), (ii) and (iv) are possible and we have
that
$$
Q(a)=a^{p+1}(1-a)^{p+1}C_1 + a^p(1-a)^p(C_2+C_4)\;.
$$
Recalling that $C_1,C_2,C_4$ are positive, each of the terms in the last equation is maximized (in $[0,1]$)
at $a=\frac{1}{2}$ which, as in the proof of Lemma \ref{l:equal}, implies that all part sizes are equal.
\qed

\section{More complete partite graphs}\label{sec:more}

\subsection{Complete partite graphs with at most one singleton part}

\begin{proof}[Proof of Theorem \ref{t:not-diverse}]
Fix $F=K_{a_1,\ldots,a_r}$ where $a_1 \ge a_2 \ge \cdots \ge a_r$ and $a_{r-1} \ge 2$.
We shall prove that there exists $t=t(F)$ such that for all $k > t$, $i_k(F)$,
and hence also $i(F)$, is attained by a complete $t$-partite graphon.

Let us first obtain a simple lower bound for $i_{r+1}(F)$ (hence, also a lower bound for $i_k(F)$ for all $k > r$). Consider $W[\frac{1}{r},\ldots,\frac{1}{r}]$ and observe that the induced density of $F$ in $W$ is at least
\begin{equation}\label{e:d}
	d := \frac{\binom{s}{a_1,\ldots,a_r}}{r^s}
\end{equation}
where $s=a_1+\cdots+a_r=|V(F)|$. So we have that $d \le i_{r+1}(F)$. 

Suppose $i_{t+1}(F)$ is attained by the complete partite graphon $W[x_1,\ldots,x_t]$
where $(x_1,\ldots,x_t) \in \Delta^{t-1}$ and $x_1 \ge x_2 \ge \cdots \ge x_t > 0$. We will prove that $t$ cannot be large. We require the following lemma.
\begin{lemma}\label{l:xr-1}
	For all sufficiently large $t$, it holds that $x_{r-1} \ge t^{-2/3}$.
\end{lemma}
\begin{proof}
	Let $t > (1/d)^3$ where $d \le i_{t+1}(F)$ is defined in \eqref{e:d}.
	Assume, by way of contradiction, that $x_{r-1} < t^{-2/3}$.
	Consider selecting independently at random $s$ points from $[0,1]$ and the event that the points induce a copy of $F$ in $W$. For this event to hold, we must have that there is an injective mapping $\sigma:[r]\rightarrow [t]$ such that precisely $a_i$ of the points fall in the part $\sigma(i)$ of $W$. But then there must be $i \le r-1$, such that $\sigma(i) \ge r-1$ and since $a_i \ge 2$
	for $i \le r-1$, part $\sigma(i)$, whose size is $x_{\sigma(i)}$, receives at least two points. 
	
	The probability that part $x_j$ receives at least two points is at most $x_j^2$.
	Summing this probability for all $r-1 \le j \le t$ and using the union bound, we have that
	the probability of inducing a copy of $F$ is at most
	$$
	\sum_{j=r-1}^{t}x_j^2 \le tx_{r-1}^2 < \frac{t}{t^{4/3}} = t^{-1/3} < d \le i_{t+1}(F)
	$$
	contradiction the assumption that $i_{t+1}(F)$ is attained by $W$.
\end{proof}

Since some parts of $F$ may have equal size, it will be more convenient to count certain ``colored'' copies of $F$ in $W$. Formalizing this, for $n \in \mathbb{N}$, let $c(n)$ be the number of vertex classes of $F$ having size precisely $n$ and let $\pi(F)=\prod_{n \in \mathbb{N}}c(n)!$.
For example, $\pi({K_{5,5,5,3,3,1}})=12$ since $c(5)=3$, $c(3)=2$, and $c(1)=1$.
Let $M$ be the set of all injective mapping from $[r]$ to $[t]$.
For $\sigma \in M$ consider the monomial
$P_\sigma=\prod_{i=1}^{r} x_{\sigma(i)}^{a_i}$ which is, up to a constant factor, the density of induced copies of $F$ in $W$ in which there are $a_i$ vertices in part $\sigma(i)$ of $W$. Then we have
\begin{equation}\label{e:before}
i_{t+1}(F) = \frac{\binom{s}{a_1,\ldots,a_r}}{\pi(F)} \sum_{\sigma \in M}\prod_{i=1}^{r} x_{\sigma(i)}^{a_i} = \frac{\binom{s}{a_1,\ldots,a_r}}{\pi(F)} \sum_{\sigma \in M} P_\sigma\;.
\end{equation}
We will show that if $t$ is sufficiently large, we can modify $(x_1,\ldots,x_t) \in \Delta^{t-1}$ to a shorter vector in $\Delta^{t-2}$, such that the sum in the right hand side of \eqref{e:before} becomes larger, which in turn, is a contradiction to the assumption that $i_{t+1}(F)$ is attained by $W[x_1,\ldots,x_t]$, whence $t$ cannot be large. We shall do this by comparing monomials before and after the modification. Consider changing $x_{t-1}$ to $x_{t-1}+x_t$, so $(x_1,\ldots,x_{t-1}+x_t) \in \Delta^{t-2}$ and let $W^*=W[x_1,\ldots,x_{t-2},x_{t-1}+x_t]$.

Let $M^*$ be the set of all injective mapping from $[r]$ to $[t-1]$.
We map each $\sigma \in M$ to $\sigma^* \in M^*$ as specified in the following cases.

(i) If $t \notin {\rm Im}(\sigma)$, then $\sigma^*=\sigma$.

(ii) If ${\rm Im}(\sigma) \cap \{t-1,t\}=\{t\}$, then let $\sigma^*$ be the same as $\sigma$ except that
$\sigma^*(\sigma^{-1}(t))=t-1$.

(iii) If ${\rm Im}(\sigma) \cap \{t-1,t\}=\{t-1,t\}$ and $\sigma(r) = t-1$,
then let $j \le r-1$ be the smallest index such that $j \notin {\rm Im}(\sigma)$. Then define
$\sigma^*$ to be the same as $\sigma$ except that:
$\sigma^*(r)=j$ and $\sigma^*(\sigma^{-1}(t))=t-1$.
For example, if $t=6$, $r=4$ and (using one-line notation) $\sigma=(1,6,3,5)$,
then we have ${\rm Im}(\sigma) =\{1,3,5,6\} \supset \{5,6\}$ so we are in this case.
Here we have $j=2$ so we have $\sigma^*(4)=2$ and $\sigma^*(2)=5$ so $\sigma^*=(1,5,3,2)$.

(iv) If ${\rm Im}(\sigma) \cap \{t-1,t\}=\{t-1,t\}$ and $\sigma(r) \neq t-1$,
then let $j \le r-1$ be the smallest index such that $j \notin {\rm Im}(\sigma)$. Then define
$\sigma^*$ to be the same as $\sigma$ except that $\sigma^*(\sigma^{-1}(t))=j$.
For example, if $t=6$, $r=4$ and (using one-line notation) $\sigma=(6,1,5,2)$,
then we have ${\rm Im}(\sigma)=\{1,2,5,6\}$, thus $j=3$ and $\sigma^{-1}(6)=1$
so $\sigma^*(1)=3$, hence $\sigma^*=(3,1,5,2)$.

Notice that our mapping from $\sigma$ to $\sigma^*$ is (obviously) not injective, yet:
Every $\sigma^*$ is the image of precisely one $\sigma$ corresponding to Case (i);
if $t-1 \in {\rm Im}(\sigma^*)$ then $\sigma^*$ is the image of precisely one possible $\sigma$ corresponding to Case (ii),  at most one possible $\sigma$ corresponding to Case (iii), and at most $r-1$ possible $\sigma$ corresponding to Case (iv).

The value of \eqref{e:before} for the modified graphon $W^*$ stays ``almost'' the same as its value for $W$. The change is that now the sum is over $M^*$, and all occurrences of $x_{t-1}$ are replaced with $x_{t-1}+x_t$.
In fact, if $\sigma^* \in M^*$ and $t-1 \in {\rm Im}(\sigma^*)$, and $P_{\sigma^*}$ before the change 
(note: $M^* \subset M$ so $\sigma^* \in M$ as well) was some $Qx_{t-1}^{a_i}$, then
$P_{\sigma^*}$ after the change is
$Q(x_t+x_{t-1})^{a_i}$ where $i={\sigma^*}^{-1}(t-1)$. 
For example, suppose $F=K_{3,2,2,1}$ (so $r=4$, $a_1=3$, $a_2=a_3=2$, $a_4=1$) and suppose $t=6$ and $\sigma^*=(2,5,4,3)$, then we have that before the change, $P_{\sigma^*}$ is
$Qx_5^2=x_2^3x_5^2x_4^2x_3$, and after the change, $P_{\sigma^*}$ is $x_2^3(x_5+x_6)^2x_4^2x_3 = Q(x_5+x_6)^2$. Further notice that we can expand $Q(x_t+x_{t-1})^{a_i}$ to at least two terms
$Qx_{t-1}^{a_i}$ (the original term) and $Qx_t^{a_i}$ and if $a_i \ge 2$, to at least one other term which is
at least $a_ix_{t-1}^{a_i-1}x_t > x_{t-1}^{a_i-1}x_t$.

In what follows we shall assume that $t \ge (r+1)^3$ and that $t > (1/d)^3$ as in Lemma \ref{l:xr-1}.
Since trivially $x_{t-1} < 1/(t-1)$ we have by Lemma \ref{l:xr-1} that $x_{r-1}/r > x_{t-1}$. 

Consider some $\sigma \in M$. We account for $P_\sigma$ before the change by comparing it to (one of the terms of) $P_\sigma^*$ after the change. We distinguish several cases according to ${\rm Im}(\sigma) \cap \{t-1,t\}$.

If ${\rm Im}(\sigma) \cap \{t-1,t\} = \emptyset$, then $\sigma^*=\sigma$ and $P(\sigma)$ before the change equals $P(\sigma^*)=P(\sigma)$ after the change, since $x_{t-1}$ and $x_t$ are not involved in $P(\sigma)$.

If ${\rm Im}(\sigma) \cap \{t-1,t\} = \{t-1\}$, then again $\sigma^*=\sigma$ and
so if $P(\sigma)=Qx_{t-1}^{a_i}$ before the change, then $P(\sigma^*)=Q(x_{t-1}+x_t)^{a_i}$ after the change
and expanding the latter, one of the terms is $Qx_{t-1}^{a_i}$, so $P(\sigma)$ is accounted for.

If ${\rm Im}(\sigma) \cap \{t-1,t\} = \{t\}$, then $\sigma^*$ is the same as $\sigma$ except that
if $\sigma(i)=t$ then now $\sigma^*(i)=t-1$. So, before the change we have
$P(\sigma)=Qx_t^{a_i}$ and $Q$ is free of $x_{t-1}$ and after the change we have
$P(\sigma^*)=Q(x_{t-1}+x_t)^{a_i}$ and after expansion, we can account for $P(\sigma)$ with the term
$Qx_t^{a_i}$ and note that we have not used this term before.

If ${\rm Im}(\sigma) \cap \{t-1,t\} = \{t-1,t\}$ and $\sigma(r) \neq t-1$, then let $\sigma(i)=t-1$
and note that this means that $i \le r-1$ so $a_i \ge 2$. Also, let $\sigma(\ell)=t$ (it is possible
that $\ell=r$). Recalling Case (iv) which is our case here, let $j$ be smallest index such
that $j \notin {\rm Im}(\sigma)$. Then before the change we have $P(\sigma)=Qx_{t-1}^{a_i}x_t^{a_\ell}$
where $Q$ is free of $x_j$. After the change, we have $P(\sigma^*)=Qx_j^{a_\ell}(x_{t-1}+x_t)^{a_i}$.
After expanding the latter, and since $a_i \ge 2$ we have a term of the form 
$a_iQx_j^{a_\ell}x_{t-1}^{a_i-1}x_t > Qx_j^{a_\ell}x_{t-1}^{a_i-1}x_t$.
But crucially, observe that
$$
\frac{Qx_j^{a_\ell}x_{t-1}^{a_i-1}x_t}{r} \ge \frac{Qx_{r-1}^{a_\ell}x_{t-1}^{a_i-1}x_t}{r}  > Qx_{t-1}^{a_i}x_t^{a_\ell} = P(\sigma) 
$$
where we have used that $x_{r-1}^{a_\ell}/r \ge (x_{r-1}/r)^{a_\ell} > x_{t-1}^{a_\ell} \ge x_{t-1}x_t^{a_\ell-1}$ recalling that $x_{r-1}/r > x_{t-1}$.
So, not only have we accounted for $P(\sigma)$, we have only used a $1/r$ portion of the term
$Qx_j^{a_\ell}x_{t-1}^{a_i-1}x_t$ and recall that we do not use this term more than $(r-1)$ times
as at most $r-1$ possible $\sigma$ corresponding to Case (iv) map to the same $\sigma^*$.
Also notice that the strict inequality in the last displayed equation shows that we are, in fact, gaining
following the change.

Finally, if ${\rm Im}(\sigma) \cap \{t-1,t\} = \{t-1,t\}$ and $\sigma(r) = t-1$, then let
$\sigma(\ell)=t$ and notice that since $\ell \neq r$ we have $\ell \le r-1$ so $a_\ell \ge 2$.
Recalling Case (iii) which is our case here, let $j$ be smallest index such
that $j \notin {\rm Im}(\sigma)$. Then before the change we have
$P(\sigma)=Qx_{t-1}^{a_r}x_t^{a_\ell}$ (possibly $a_r=1$) where $Q$ is free of $x_j$.
After the change, we have $P(\sigma^*)=Qx_j^{a_r}(x_{t-1}+x_t)^{a_\ell}$.
After expanding the latter, and since $a_\ell \ge 2$ we have a term of the form 
$a_\ell Qx_j^{a_r}x_{t-1}^{a_\ell-1}x_t > Qx_j^{a_r}x_{t-1}^{a_\ell-1}x_t$.
But crucially, observe that
$$
\frac{Qx_j^{a_r}x_{t-1}^{a_\ell-1}x_t}{r} \ge \frac{Qx_{r-1}^{a_r}x_{t-1}^{a_\ell-1}x_t}{r}  > Qx_{t-1}^{a_r}x_t^{a_\ell} = P(\sigma) 
$$
where we have used again that $x_{r-1}^{a_r}/r \ge (x_{r-1}/r)^{a_r} > x_{t-1}^{a_r}$
and that $x_{t-1} \ge x_t$.
So, not only have we accounted for $P(\sigma)$, we have only used a $1/r$ portion of the term
$Qx_j^{a_r}x_{t-1}^{a_\ell-1}x_t$ and recall that we do not use this term more than once,
at most one possible $\sigma$ corresponding to Case (iii) maps to a particular $\sigma^*$.
(Nevertheless, recall from the previous case that in addition we might have at most $r-1$ possible $\sigma$ corresponding to Case (iv) that map to $\sigma^*$ which is why we have used an $r$-fraction of that term.)
\end{proof}

\subsection{Strongly unbalanced complete partite graphs}

\begin{proof}[Proof of Theorem \ref{t:strongly-unbalanced}]
	Suppose that $F=K_{a_1,\ldots,a_r}$ is strongly unbalanced.
	We will prove Theorem \ref{t:strongly-unbalanced} for the case of $i(F)$; the case of $i_k(F)$ is analogous and simpler. Since $F$ is strongly unbalanced, we have that $(a_i-a_j)^2 > (a_i+a_j)$ for $1 \le i < j \le r$.
	In particular, assuming $a_i \ge a_{i+1}$ for $1 \le i \le r-1$, we have that all part sizes of $F$ are distinct and in particular $a_{r-1} \ge 2$, so by Theorem \ref{t:not-diverse}, there is $t=t(F)$ such that for all $k > t$,
	$i_k(F)$ is attained by a complete $t$-partite graphon $W$, and hence $i(F)$ is also obtained by $W$.
	
	Let $W=W(x_1,\ldots,x_t)$ where $(x_1,\ldots,x_t) \in \Delta^{t-1}$ and $x_i > 0$ for $1 \le i \le t$.
	In particular, there is $\gamma=\gamma(F)$ such that $x_i \ge \gamma$ for $1 \le i \le t$.
	We must prove that $x_i \neq x_j$ for all $1 \le i < j \le t$.
	
	Recall that the induced density of $F$ in $W$ is given by a symmetric polynomial in $x_1,\ldots,x_t$ which
	we shall denote here by $P_F$. Assume, by way of contradiction, that
	$x_i = x_j$ for some $i < j$. We shall derive a contradiction by showing that
	for some small $\delta > 0$,
	$$
	P^*=P(x_1,\ldots,x_i+\delta,x_{i+1},\ldots,x_{j-1},x_i-\delta,x_{j+1},\ldots,x_t) > P(x_1,\ldots,x_t)\;.
	$$
	
	Let $\sigma: [r] \rightarrow [t]$ be injective and consider the corresponding monomial in $P$, which, recall by the proof of Theorem \ref{t:not-diverse}, is $P_\sigma=\prod_{i=1}^{r} x_{\sigma(i)}^{a_i}$ \footnote{Recall that each monomial is multiplied by the same constant $\binom{a_1+\cdots+a_r}{a_1,\ldots,a_r}$; in our case we have $\pi(F)=1$ since all part sizes of $F$ are distinct.}.
	Each $P(\sigma)$ is a monomial in $x_1,\ldots,x_t$, where only the variables $x_\ell$ for
	$\ell \in {\rm Im}(\sigma)$ participate.
	We shall show that if we replace $x_i$ with $x_i+\delta$ and replace $x_j=x_i$ with $x_i-\delta$
	then either $P(\sigma)$ does not change, or there is a pair $\sigma,\sigma'$ such that $P(\sigma)+P(\sigma')$ increases, hence we shall arrive at a contradiction.
	
	Consider first the case that ${\rm Im}(\sigma) \cap \{i,j\} = \emptyset$. In this case,
	$P(\sigma)$ does not change after the change of variables.
	
	Consider next the case that ${\rm Im}(\sigma) \cap \{i,j\} = \{i\}$ or ${\rm Im}(\sigma) \cap \{i,j\} = \{j\}$.
	Note that if ${\rm Im}(\sigma) \cap \{i,j\} = \{i\}$ and $\sigma^{-1}(i)=\ell$, then $\sigma^*$ which is defined
	the same as $\sigma$ except that $\sigma^*(\ell)=j$ has that $j \in {\rm Im}(\sigma^*)$,
	so there is a matching between these two types of injections (those whose image contains $i$ and not $j$ those whose image contains $j$ and not $i$).
	Let us observe $P(\sigma)+P(\sigma^*)$ before and after the change of variables.
	Before the change, this sum was of the form $Qx_i^{a_\ell}+Qx_j^{a_\ell} = 2Qx_i^{a_\ell}$ (recall, $x_i=x_j$)
	where $Q$ is free of $x_i$ and $x_j$. After the change, it is of the form
	$Q(x_i+\delta)^{a_\ell}+Q(x_i-\delta)^{a_\ell}$. But notice that
	$$
	2x_i^{a_\ell} < (x_i+\delta)^{a_\ell}+(x_i-\delta)^{a_\ell}
	$$
	so there is an increase, as desired.
	
	Finally, consider the case ${\rm Im}(\sigma) \cap \{i,j\} = \{i,j\}$ and let $\sigma^*$
	be defined the same as $\sigma$ except that $\sigma^*(\sigma^{-1}(i))=j$,
	$\sigma^*(\sigma^{-1}(j))=i$. Let $\sigma^{-i}(i)=\ell$ and $\sigma^{-1}(j)=m$.
	Let us observe $P(\sigma)+P(\sigma^*)$ before and after the change of variables.
	Before the change, this sum was of the form $2Qx_i^{a_\ell+a_m}$
	where $Q$ is free of $x_i$ and $x_j$.
	After the change, it is of the form $Q(x_i+\delta)^{a_\ell}(x_i-\delta)^{a_m} +
	Q(x_i-\delta)^{a_\ell}(x_i+\delta)^{a_m}$. It suffices to prove that for small enough $\delta$,
	$$
	2x_i^{a_\ell+a_m} < (x_i+\delta)^{a_\ell}(x_i-\delta)^{a_m} + (x_i-\delta)^{a_\ell}(x_i+\delta)^{a_m}\;.
	$$
	Indeed,
	\begin{align*}
	& (x_i+\delta)^{a_\ell}(x_i-\delta)^{a_m} + (x_i-\delta)^{a_\ell}(x_i+\delta)^{a_m} - 2x_i^{a_\ell+a_m}\\
	= ~& \delta^2[a_m(a_m-1)+a_\ell(a_\ell-1)-2a_m a_\ell]x_i^{a_m+a_\ell-2} + \delta^4Q^*(x_i,\delta)
	\end{align*}
	where $Q^*$ is a polynomial in $\delta$ and $x_i$ (which may evaluate to negative).
	Hence, it suffices to prove that  $a_m(a_m-1)+a_\ell(a_\ell-1)-2a_m a_\ell > 0$.
	But this translates to showing that $(a_\ell - a_m)^2 > a_\ell + a_m$ which indeed holds since $F$ is
	strongly unbalanced.	
\end{proof}

\subsection{Proof of Proposition \ref{prop:smallest}}

Recall that we aim to show the existence of a complete partite graph $F$ for which $i(F)$ is rational,
yet $i_k(F)$ is irrational for some $k$, and determine the smallest such graph.
Let us first rule out some candidates, so that we are left with a smallest possible candidate, and then show that it does, in fact, satisfy the constraints.

Any complete bipartite graph is ruled out, as we have mentioned in the introduction that for such graphs $F$, Brown and Sidorenko \cite{BS-1994} proved that $i(F)=i_3(F)$ (and trivially $i_2(F)=0$).
Also, Theorem \ref{t:turan} rules out any Tur\'an graph $T(s,r)$ with $s \le 3r+1$, which already rules out
all such graphs on at most $10$ vertices (since $r \ge 3$).
The smallest non-bipartite complete partite graph that is not a Tur\'an graph has $5$ vertices and is $K_{3,1,1}$. Recall from the introduction that $i(K_{3,1,1})=\frac{216}{625}$ 
as proved in \cite{LPSS-2023} (we note that the lower bound construction and a presumably tight
flag algebra upper bound were already given in \cite{EL-2015}). So, it suffices to prove that
$i_k(K_{3,1,1})$ is irrational for some $k$ (and obviously $k \ge 4$).

Let us first consider the case $k=4$ and the polynomial exhibiting the induced density of $K_{3,1,1}$ in the graphon $W[x,y,z]$ which is
$$
20(x^3yz+xy^3z+xyz^3) = 20xyz(x^2+y^2+z^2)\;.
$$
It is easily shown using Lagrange multipliers to maximize for $(x,y,z) = (\frac{1}{3},\frac{1}{3},\frac{1}{3})$
at which point we have $i_4(K_{3,1,1})= \frac{20}{81}$, still rational.

Consider next the case $k=5$ where the polynomial in terms of $(x,y,z,w)$ is
$$
20(x^3yz+xy^3z+xyz^3+x^3yw+xy^3w+xyw^3+x^3wz+xw^3z+xwz^3+w^3yz+wy^3z+wyz^3)
$$
which is still solvable analytically using Lagrange multipliers, though quite tedious to perform manually.
The following Maple worksheet provides the complete solution \url{https://github.com/raphaelyuster/turan-inducibility/blob/main/K311.mw}. A global maximum is
at $(\alpha,\alpha,\alpha,1-3\alpha)$  (and, by symmetry, there are four global maxima) where
$\alpha = (19- \sqrt{41})/80=0.1574...$. Evaluated at these points we obtain
$$
i_5(K_{3,1,1}) = \frac{329847+ 1107\sqrt{41}}{1280000} = 0.2632...
$$
Note that the diagonal point $(\frac{1}{4},\frac{1}{4},\frac{1}{4},\frac{1}{4})$ only evaluates to $\frac{15}{64}$.
\qed


\begin{thebibliography}{10}
	
	\bibitem{AHKT-2020}
	N.~Alon, D.~Hefetz, M.~Krivelevich, and M.~Tyomkyn.
	\newblock Edge-statistics on large graphs.
	\newblock {\em Combinatorics, Probability and Computing}, 29(2):163--189, 2020.
	
	\bibitem{AS-2016}
	N.~Alon and C.~Shikhelman.
	\newblock Many {$T$} copies in {$H$}-free graphs.
	\newblock {\em Journal of Combinatorial Theory, Series B}, 121:146--172, 2016.
	
	\bibitem{BHLP-2016}
	J.~Balogh, P.~Hu, B.~Lidick{\`y}, and F.~Pfender.
	\newblock Maximum density of induced 5-cycle is achieved by an iterated blow-up
	of 5-cycle.
	\newblock {\em European Journal of Combinatorics}, 52:47--58, 2016.
	
	\bibitem{BEHJ-1995}
	B.~Bollob{\'a}s, Y.~Egawa, A.~Harris, and G.~Jin.
	\newblock The maximal number of induced $r$-partite subgraphs.
	\newblock {\em Graphs and Combinatorics}, 11(1):1--19, 1995.
	
	\bibitem{BS-1994}
	J.~I. Brown and A.~Sidorenko.
	\newblock The inducibility of complete bipartite graphs.
	\newblock {\em Journal of Graph Theory}, 18(6):629--645, 1994.
	
	\bibitem{DHMNS-2013}
	S.~Das, H.~Huang, J.~Ma, H.~Naves, and B.~Sudakov.
	\newblock A problem of {E}rd\"os on the minimum number of $k$-cliques.
	\newblock {\em Journal of Combinatorial Theory, Series B}, 103(3):344--373,
	2013.
	
	\bibitem{EFR-1986}
	P.~Erd\H{o}s, P.~Frankl, and V.~R{\"o}dl.
	\newblock The asymptotic number of graphs not containing a fixed subgraph and a
	problem for hypergraphs having no exponent.
	\newblock {\em Graphs and Combinatorics}, 2(1):113--121, 1986.
	
	\bibitem{erdos-1962}
	P.~Erd{\H{o}}s.
	\newblock On the number of complete subgraphs contained in certain graphs.
	\newblock {\em Magyar Tudom{\'a}nyos Akad{\'e}mia Matematikai Kutat{\'o}
		Int{\'e}zet{\'e}nek K{\"o}zlem{\'e}nyei}, 7(3):459--464, 1962.
	
	\bibitem{EL-2015}
	C.~Even-Zohar and N.~Linial.
	\newblock A note on the inducibility of $4$-vertex graphs.
	\newblock {\em Graphs and Combinatorics}, 31(5):1367--1380, 2015.
	
	\bibitem{exoo-1986}
	G.~Exoo.
	\newblock Dense packings of induced subgraphs.
	\newblock {\em Ars Combinatoria}, 22:5--10, 1986.
	
	\bibitem{FHL-2017}
	J.~Fox, H.~Huang, and C.~Lee.
	\newblock A solution to the inducibility problem for almost all graphs.
	\newblock {\em preprint}, 2017.
	
	\bibitem{FS-2020}
	J.~Fox and L.~Sauermann.
	\newblock A completion of the proof of the edge-statistics conjecture.
	\newblock {\em Advances in Combinatorics}, page~4, 2020.
	
	\bibitem{goodman-1959}
	A.~W. Goodman.
	\newblock On sets of acquaintances and strangers at any party.
	\newblock {\em The American Mathematical Monthly}, 66(9):778--783, 1959.
	
	\bibitem{grzesik-2012}
	A.~Grzesik.
	\newblock On the maximum number of five-cycles in a triangle-free graph.
	\newblock {\em Journal of Combinatorial Theory, Series B}, 102(5):1061--1066,
	2012.
	
	\bibitem{HHN-2014}
	H.~Hatami, J.~Hirst, and S.~Norine.
	\newblock The inducibility of blow-up graphs.
	\newblock {\em Journal of Combinatorial Theory, Series B}, 109:196--212, 2014.
	
	\bibitem{HHKNR-2013}
	H.~Hatami, J.~Hladk{\`y}, D.~Kr{\'a}l, S.~Norine, and A.~Razborov.
	\newblock On the number of pentagons in triangle-free graphs.
	\newblock {\em Journal of Combinatorial Theory, Series A}, 120(3):722--732,
	2013.
	
	\bibitem{HT-2018}
	D.~Hefetz and M.~Tyomkyn.
	\newblock On the inducibility of cycles.
	\newblock {\em Journal of Combinatorial Theory, Series B}, 133:243--258, 2018.
	
	\bibitem{hirst-2014}
	J.~Hirst.
	\newblock The inducibility of graphs on four vertices.
	\newblock {\em Journal of Graph Theory}, 75(3):231--243, 2014.
	
	\bibitem{KNV-2019}
	D.~Kr{\'a}l, S.~Norin, and J.~Volec.
	\newblock A bound on the inducibility of cycles.
	\newblock {\em J. Combin. Theory Ser. A}, 161:359--363, 2019.
	
	\bibitem{KST-2019}
	M.~Kwan, B.~Sudakov, and T.~Tran.
	\newblock Anticoncentration for subgraph statistics.
	\newblock {\em Journal of the London Mathematical Society}, 99(3):757--777,
	2019.
	
	\bibitem{LP-2018}
	B.~Lidick{\`y} and F.~Pfender.
	\newblock Pentagons in triangle-free graphs.
	\newblock {\em European Journal of Combinatorics}, 74:85--89, 2018.
	
	\bibitem{LPSS-2023}
	H.~Liu, O.~Pikhurko, M.~Sharifzadeh, and K.~Staden.
	\newblock Stability from graph symmetrisation arguments with applications to
	inducibility.
	\newblock {\em Journal of the London Mathematical Society}, 108(3):1121--1162,
	2023.
	
	\bibitem{LMR-2023}
	X.~Liu, D.~Mubayi, and C.~Reiher.
	\newblock The feasible region of induced graphs.
	\newblock {\em Journal of Combinatorial Theory, Series B}, 158:105--135, 2023.
	
	\bibitem{lovasz-2012}
	L.~Lov{\'a}sz.
	\newblock {\em Large networks and graph limits}, volume~60.
	\newblock American Mathematical Soc., 2012.
	
	\bibitem{MMNT-2019}
	A.~Martinsson, F.~Mousset, A.~Noever, and M.~Truji{\'c}.
	\newblock The edge-statistics conjecture for $\ell \ll k^{6/5}$.
	\newblock {\em Israel Journal of Mathematics}, 234(2):677--690, 2019.
	
	\bibitem{PST-2019}
	O.~Pikhurko, J.~Slia{\v{c}}an, and K.~Tyros.
	\newblock Strong forms of stability from flag algebra calculations.
	\newblock {\em Journal of Combinatorial Theory, Series B}, 135:129--178, 2019.
	
	\bibitem{PV-2013}
	O.~Pikhurko and E.~R. Vaughan.
	\newblock Minimum number of $k$-cliques in graphs with bounded independence
	number.
	\newblock {\em Combinatorics, Probability and Computing}, 22(6):910--934, 2013.
	
	\bibitem{PG-1975}
	N.~Pippenger and M.~C. Golumbic.
	\newblock The inducibility of graphs.
	\newblock {\em Journal of Combinatorial Theory, Series B}, 19(3):189--203,
	1975.
	
	\bibitem{szemeredi-1978}
	E.~Szemer\'edi.
	\newblock Regular partitions of graphs.
	\newblock In {\em Probl\`emes combinatoires et th\'eorie des graphes ({C}olloq.
		{I}nternat. {CNRS}, {U}niv. {O}rsay, {O}rsay, 1976)}, volume 260 of {\em
		Colloq. Internat. CNRS}, pages 399--401. CNRS, Paris, 1978.
	
	\bibitem{turan-1941}
	P.~Tur{\'a}n.
	\newblock On an extremal problem in graph theory.
	\newblock {\em Mat. Fiz. Lapok}, 48:436--452, 1941.
	
	\bibitem{ueltzen-2024}
	R.~Ueltzen.
	\newblock Characterizing graphs with high inducibility.
	\newblock {\em arXiv preprint arXiv:2411.17362}, 2024.
	
	\bibitem{flagmatic-site}
	E.~Vaughan.
	\newblock Flagmatic.
	\newblock \url{https://github.com/emil79/flagmatic}, 2012.
	
	\bibitem{yuster-2019}
	R.~Yuster.
	\newblock On the exact maximum induced density of almost all graphs and their
	inducibility.
	\newblock {\em Journal of Combinatorial Theory, Series B}, 136:81--109, 2019.
	
	\bibitem{zykov-1949}
	A.~A. Zykov.
	\newblock On some properties of linear complexes.
	\newblock {\em Matematicheskii Sbornik}, 66(2):163--188, 1949.
	
\end{thebibliography}
\end{document}